\documentclass[11pt]{amsart}
\usepackage{xcolor}
\usepackage{graphicx}
\usepackage{dsfont}
\usepackage{amsmath}
\usepackage[hidelinks]{hyperref}
\usepackage{float}
\usepackage{bbm}
\usepackage{enumerate}
\usepackage{todonotes}
\usepackage{mathtools}
\mathtoolsset{showonlyrefs}

\newtheorem{theorem}{Theorem}[section]
\newtheorem{lemma}{Lemma}[section]
\newtheorem{prop}{Proposition}[section]
\newtheorem{definition}{Definition}[section]

\newtheorem{remark}{Remark}[section]

\newcommand{\be}{\begin{equation}}
\newcommand{\ee}{\end{equation}}
\newcommand{\bea}{\begin{eqnarray}}
\newcommand{\eea}{\end{eqnarray}}
\newcommand{\barr}{\begin{array}}
\newcommand{\earr}{\end{array}}

\DeclareMathOperator{\sgn}{sgn}

\newcommand{\bpar}{\begin{equation} \left\{ \begin{array}{lll}}
\newcommand{\epar}{ \end{array}\right. \end{equation} }
\newcommand{\eparn}{ \end{array} \right.}

\newcommand\leftmat{\left(\begin{array}{cc}}
\newcommand\rightmat{\end{array}\right)}
\newcommand\leftvec{\left(\begin{array}{c}}
\newcommand\rightvec{\end{array}\right)}
\newcommand\re{{\rm e}}
\newcommand\la{{\lambda}}
\newcommand\R{{\mathbb R}}
\newcommand\C{{\mathbb C}}
\newcommand\Z{{\mathbb Z}}
\newcommand\N{{\mathbb N}}

\def\Xint#1{\mathchoice
   {\XXint\displaystyle\textstyle{#1}}%
   {\XXint\textstyle\scriptstyle{#1}}%
   {\XXint\scriptstyle\scriptscriptstyle{#1}}%
   {\XXint\scriptscriptstyle\scriptscriptstyle{#1}}%
   \!\int}
\def\XXint#1#2#3{{\setbox0=\hbox{$#1{#2#3}{\int}$}
     \vcenter{\hbox{$#2#3$}}\kern-.5\wd0}}

\def\dashint{\Xint-}

\providecommand{\bigohsymbolonly}{\operatorname{O}}
\providecommand{\bigoh[1]}{\bigohsymbolonly\!\left(#1\right)}

\providecommand{\Lebesgue}{\operatorname{L}} 
\providecommand{\ContinuousSpace}{\operatorname{C}} 

\providecommand{\Ltwozo}{\Lebesgue^2(0,1)}       
\providecommand{\Czo}{\ContinuousSpace[0,1]}     
\providecommand{\Czopen}{\ContinuousSpace(0,1)}

\providecommand{\argdot}{{}\cdot{}}

\renewcommand{\Re}{\operatorname{Re}}

\title[Jumps and cusps: a new revival effect]{Jumps and cusps: a new revival effect in local dispersive PDEs}
\author[Boulton, Farmakis, Pelloni, Smith]{
L. Boulton$^\beta$, G. Farmakis$^\varphi$,  B. Pelloni$^\pi$ \& D.A. Smith$^\sigma$}

\address{
$^{\beta,\varphi,\pi}$Heriot-Watt University \& Maxwell Institute  for the Mathematical Sciences, Edinburgh, Scotland. $^\sigma$Yale-NUS College \& National University of Singapore, Singapore.
}

\email{$^\beta$L.Boulton@hw.ac.uk, $^\varphi$G.Farmakis@hw.ac.uk, $^\pi$B.Pelloni@hw.ac.uk, $^\sigma$dave.smith@yale-nus.edu.sg}

\date{Version: 1/3/24}

\begin{document}
\maketitle

\begin{abstract}
We study the presence of a non-trivial revival effect in the solution of linear dispersive boundary value problems for two benchmark models which arise in applications: the Airy equation and the dislocated Laplacian Schr{\"o}dinger equation. In both cases, we consider boundary conditions of Dirichlet-type. We prove that, at suitable times, jump discontinuities in the initial profile are revived in the solution, not only as jump discontinuities but also as logarithmic cusp singularities. We explicitly describe these singularities and show that their formation is due to interactions between the symmetries of the underlying spatial operators with the periodic Hilbert transform.
\end{abstract}


\pagebreak

\section{Introduction}
The phenomenon of {\em revivals} in linear dispersive periodic boundary value problems, known in the literature also as {\em Talbot effect}, has been studied for over 30 years and is now well understood. The currently accepted, heuristic definition of this phenomenon states that, at specific values of the time variable called {\em rational times}, the solution is  a linear combination of a finite number of translated copies of the initial condition.
This behaviour is brought into sharp relief when this initial condition has jump discontinuities, as the jumps then propagate through the solution at every rational time. 
For a comprehensive introduction to this phenomenon, with historical and bibliographical notes,  see \cite{erdougan2016dispersive}. 
We will refer to this type of phenomenon
  as the {\bf {\em periodic revival} }property.

Recently, a novel manifestation of revivals has been observed in the context of linear {\em nonlocal} equations, which notably includes the linearised Benjamin-Ono equation, \cite{boulton2020new}. In this manifestation, the initial jump discontinuities give rise not only to  jumps but also to logarithmic cusp singularities appearing in the solution. We will call this the {\bf {\em cusp revival}} property.

In parallel to the discovery of this novel effect, a relaxation of the notion of revivals has also been considered recently, \cite{BFP, boulton2023phenomenon,GFthesis}. This relaxation enables the study of a larger class of problems and it is clearly important for the rigorous mathematical description of revivals from a more general point of view. We will refer to this as the {\bf {\em weak revival}} property, of either periodic or cusp type.

\smallskip
In the present paper we show that  {\em local} linear dispersive equations can exhibit the weak cusp revival property, when non-periodic boundary conditions are imposed. The fact that the boundary conditions are not periodic is crucial to the formation of this more complex type of revival structure. 

We consider two canonical models which are naturally linked to applications. The first model is the {\em Airy equation\footnote{ The name {\em Airy equation} is used on account of the fact that the fundamental solution on $\R$ of the given PDE is the Airy function $\operatorname{Ai}(x/\sqrt[3]{3t})$.} with Dirichlet-type boundary conditions} on $[0,1]$. The boundary value problem is given by
\begin{equation} \label{pseD1} \tag{A}
\begin{aligned}
&u_t(x,t)+u_{xxx}(x,t)=0, &  \qquad x\in(0,1), \;t\in \mathbb{R}, \\
&u(0,t)=0,\quad u(1,t)=0,\quad u_x(1,t)= u_x(0,t), &\qquad t\in \mathbb{R},\nonumber\\
&u(x,0)=u_0(x), &\qquad  x\in(0,1).
\end{aligned}
\end{equation}
An explicit formula for the solution of \eqref{pseD1} is important for accurate  numerical simulations of  wave motion, see the discussion in \cite{bsz}. 

The second model is the {\em linear Schr\"odinger equation for the dislocated Laplacian with Dirichlet boundary conditions} on $[0,1]$, given by
\begin{equation} \label{disPDE} \tag{D}
\begin{aligned}
&u_t(x,t)+iu_{xx}(x,t)=0, && x\in(0,b),\;\;t\in \mathbb{R},
\\
&u_t(x,t)-iu_{xx}(x,t)=0, && x\in(b,1),\;\;t\in \mathbb{R},
\\
&u(0,t)=u(1,t)=0,
\qquad&& t\in \mathbb{R},\\
&u(b^+,t)=u(b^-,t),\quad u_x(b^+,t)=-u_x(b^-,t), && t\in \mathbb{R}, \\
&u(x,0)=u_0(x),  && x\in(0,1).
\end{aligned}
\end{equation}
Here, $b\in(0,1)$ is the dislocation. The boundary value problem \eqref{disPDE} is a prototype for the study of cloaking in meta-materials, see \cite{JussiDavid2018} and references therein.

\smallskip

Our main contributions are stated in Section~\ref{section2}. In a nutshell, we establish that if $u_0$ is a piecewise Lipschitz function,
then the solution at appropriate rational times is a superposition of the following three components:
\begin{enumerate}[a)]
\item \label{pra}
a finite linear combination of translated copies of an auxiliary function $G_{u_0}$, reviving any jump of $u_0$;
 \item \label{prb}
a finite linear combination of copies of the periodic Hilbert transform of $G_{u_0}$, producing a finite number of logarithmic cusps corresponding to the jumps of $u_0$;
 \item \label{prc}
a continuous function of the spatial variable $x$.
 \end{enumerate}
The function $G_{u_0}$  is given explicitly in terms of the initial condition $u_0$. 

In the case of the boundary value problem \eqref{disPDE}, this representation holds separately in the sub-intervals $(0,b)$ and $(b,1)$.  Indeed,  the dislocation at $x=b$ acts as a barrier, so that any discontinuity initially in one of the sub-intervals gives rise to jumps and cusps only there.

\begin{figure}[ht]
  \centering
    \includegraphics[width=.8\linewidth]{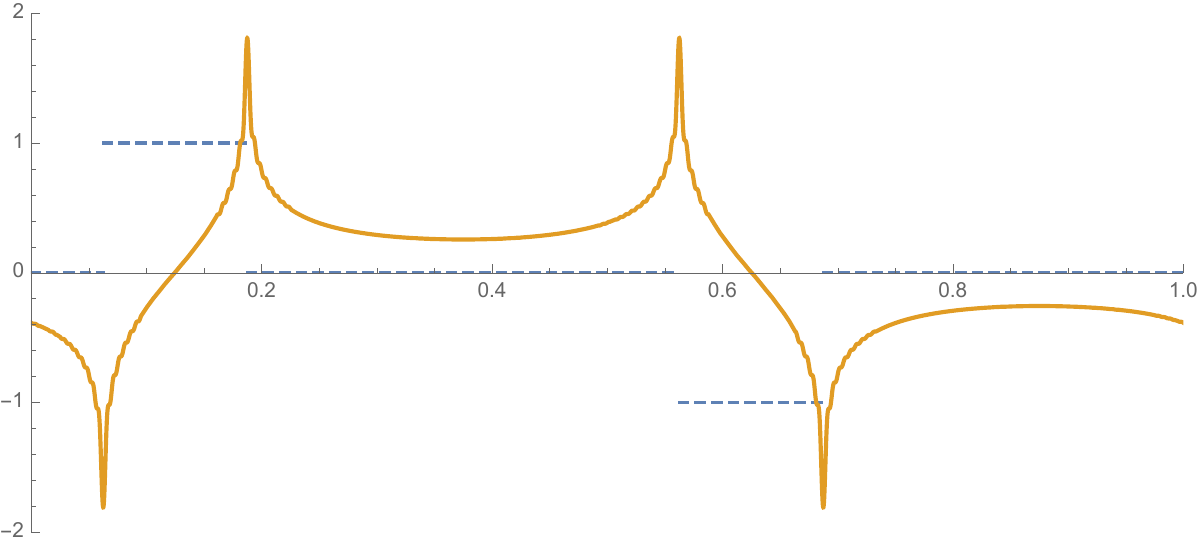}
  \caption{A step function (dashed) and its periodic Hilbert transform (solid). When a given function has an isolated jump discontinuity, its Hilbert transform displays a cusp with a singularity of logarithmic order, see \eqref{logsing}. (The cusps actually extend to infinity in either direction).}
    \label{fig0}
   \end{figure}

The cusps observed in the solutions arise from the singularities of the periodic Hilbert transform of $G_{u_0}$. Indeed, the periodic Hilbert transform of any function with isolated jump discontinuities has logarithmic cusps, such as those shown in Figure~\ref{fig0}.

The boundary value problems \eqref{pseD1} and \eqref{disPDE} have the following characteristics in common:
\begin{enumerate}[i)]
\item \label{pr1}  the underlying linear operators are determined by boundary conditions which cannot be re-formulated in terms of periodic ones;
\item\label{pr2} the spatial linear operators are indefinite, in the sense that their spectrum  is discrete and accumulates at $\pm\infty$;
\item\label{pr3}  the solution can be formulated in terms of the periodic Hilbert transform of a function explicitly determined by the initial condition.
\end{enumerate}
On account of \ref{pr1}), it is surprising that either problem exhibits any form of revivals. This is also surprising given the property \ref{pr2}), since a change of sign in just one eigenvalue can result in a mistuning of the  revival property, \cite{BFP}.  It is even more surprising that, as we show 
below,  there is a common technique to derive the property \ref{pr3}) for both equations.

We highlight an important distinction between the case of the linear Benjamin-Ono equation, studied  in \cite{boulton2020new}, and the ones that we consider here. In the former, the cusps naturally appear as a consequence of the presence of the Hilbert transform in the formulation of the PDE. Such nonlocal term is not present in the equations considered in this paper. As we show below, the cusp revival property arises as a consequence of the interaction of the symmetries of the leading part of the spatial operators with those of the periodic Hilbert transform.


\section{Statement of results} \label{section2}

\subsection*{The Dirichlet-type Airy equation} 
Our main contribution about the boundary value problem \eqref{pseD1} is summarised in the next theorem, which implies that the solution exhibits weak cusp revivals. Specifically, up to a continuous perturbation, for all $t\in\mathbb{R}$ the solution is represented by an infinite series in terms of the explicit exponential functions ${\displaystyle\re^{i(2n-\frac13)\pi x}}$ with temporal evolution of the form ${\displaystyle\re^{i(2n-\frac13)^3\pi^3 t}}$. At the specific rational times \({\displaystyle t=\frac p  {\pi^2 q}}\), the infinite superposition of these components can be rearranged into a finite summation of functions which have simple explicit expressions. This finite summation can be further decomposed into terms explicitly describing  the jumps and cusps characteristic of the  cusp revival property.

\begin{theorem}\label{lkdvthm2}
 Let $u_0:[0,1]\longrightarrow \R$ be piecewise Lipschitz.
    Then, for all $t\in\R$, the solution of the boundary value problem~\eqref{pseD1} admits the representation
    $$
        u(x,t)=U_{\mathcal R}(x,t)+U_{\mathcal C}(x,t)
    $$
    where
    $U_{\mathcal C}(\argdot,t)\in\Czopen$
     and
    \begin{equation} \label{Ldef0}
        U_{\mathcal R}(x,t)=\sum_{n=1}^\infty 2 \operatorname{Re} \left[\widetilde{u_0}(n)
        \re^{i(2n-\frac 1 3 )^3\pi^3t}\re^{i (2n-\frac 1 3 )\pi x} \right],
    \end{equation}
    with
    \be \label{eqn:defn.u0hat}
   \widetilde{u_0}(n)=    \int_0^1 \left[u_0(y)+u_0(1)+u_0(0)\right] \re^{-i(2n-\frac 1 3 )\pi y}\mathrm{d}y.
    \ee
If $p,q\in\N$ are coprime, then 
    \begin{equation} \label{finrevf2NEW}
        U_{\mathcal R}\left(x,\frac p {q \pi^2}\right) = \Re\left\{ \re^{-i\frac{\pi}{27}\left(9x+\frac pq\right)} [L_1(x)+iL_2(x)-L_3] \right\}
    \end{equation}
    where
    \begin{align*}
        L_1(x)
        &= \sum_{k=0}^{q-1} d_k^{p,q}
        G_{u_0}\left(x+\frac p{3q}-\frac{k}{q}\right),
        \\
        L_2(x)
        &= \sum_{k=0}^{q-1} d_k^{p,q}
        [{\mathcal H}G_{u_0}]\left(x+\frac p{3q}-\frac{k}{q}\right), \\
        L_3
        &=
            \int_0^1G_{u_0}(y)\mathrm{d}y 
    \end{align*}
and
    \begin{equation} \label{v0NEW}
        G_{u_0}(x)=[u_0(x)+u_0(1)+u_0(0)]\re^{i\frac{\pi x}3}.
    \end{equation}    
    Here $\mathcal{H}$ denotes the $1$-periodic Hilbert transform  and $d_k^{p,q}$ are scalar coefficients independent of $u_0$.
\end{theorem}

Note that, if the initial condition $u_0$ is compatible with the homogeneous Dirichlet boundary conditions in \eqref{pseD1} so that $u_0(0)=u_0(1)=0$, then $ \widetilde{u_0}(n)$ is simply the generalised Fourier coefficient of the initial condition with respect to the orthogonal sequence $\{\re^{i(2n-\frac13)\pi x}\}_{n\in\mathbb{N}}$.

The function $G_{u_0}$ in \eqref{v0NEW} has jump discontinuities if and only if $u_0$ does. Its periodic Hilbert transform, appearing in the term $L_2$, will display a logarithmic cusp whenever $u_0$ has a jump discontinuity. These will then combine with the jump discontinuities of the term $L_1$ ($L_3$ is constant), to form the solution profile characteristic of the cusp revival phenomenon. 

\begin{figure}
        \centering
        \includegraphics[width=0.98\linewidth]{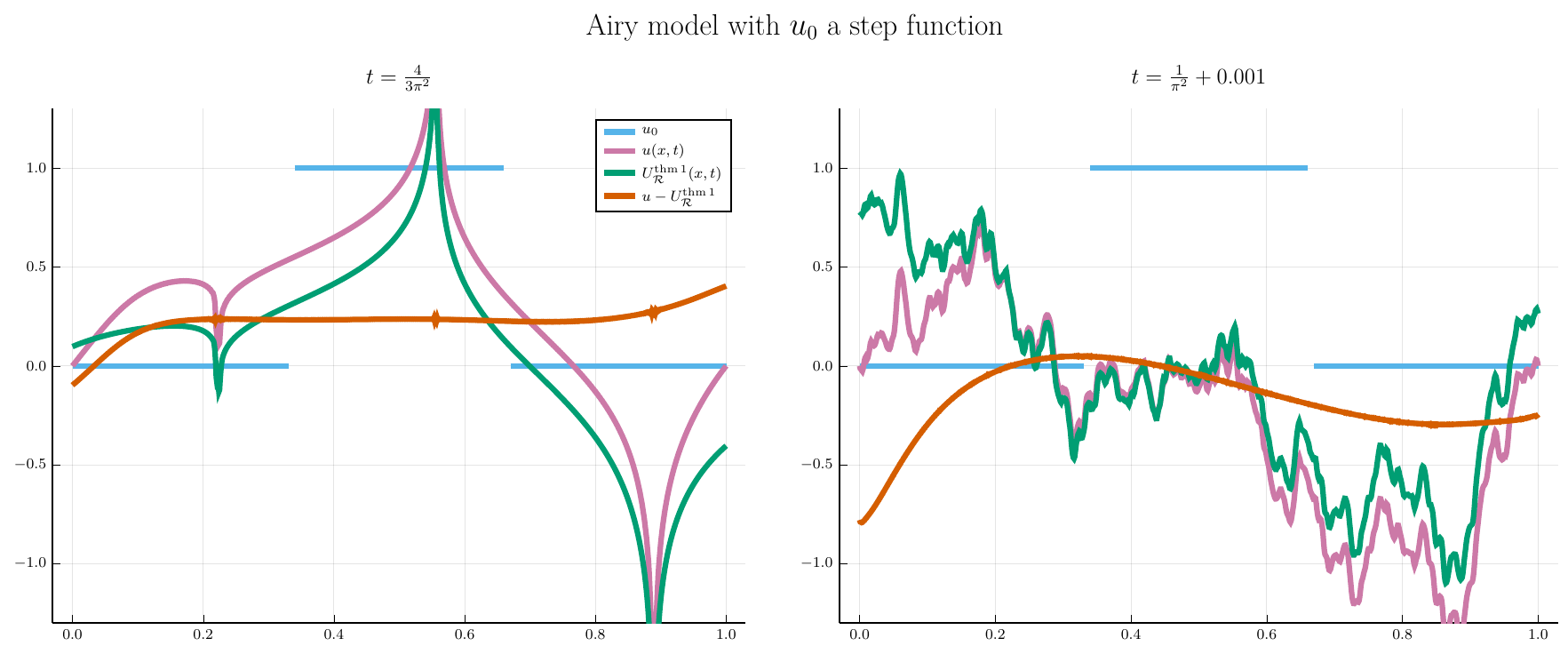}
        \caption{
For $u_0$ a step function (blue) the graphs show the solution of \eqref{pseD1} (purple), the function $U_{\mathcal R}$ (green) and the continuous function $U_{\mathcal C}$ (brick). See Theorem~\ref{lkdvthm2}. The horizintal axis is $x\in[0,1]$. The time $t$ is fixed rational on the left and fixed irrational on the right. These graphs were constructed from an expansion in the first 600 eigenfunctions, approximated numerically.
       }
        \label{fig:intro.A}
    \end{figure}

In Figure~\ref{fig:intro.A}, we show a numerical approximation of the solution corresponding to an initial $u_0$ given by a step function. We also shown both $U_{\mathcal R}$ and $U_{\mathcal C}$.  The graph on the left hand side illustrates the appearance of  cusps in both the solution and $U_{\mathcal R}$, at a time that is a rational multiple of $\pi^{-2}$ matching the expression \eqref{finrevf2NEW}. The graph on the right hand side is indicative of a possible fractalisation phenomenon at irrational times. 

The proof of Theorem~\ref{lkdvthm2}, given in Section~\ref{sec4}, is not a consequence of the revival property for periodic boundary conditions, as would be the case for second order problems, \cite{farmakis2024}. Indeed, while the boundary value problem for the Airy equation with periodic boundary conditions exhibits periodic revivals, \cite{olver2010dispersive}, the boundary value problem \eqref{pseD1} is not trivially related to it, despite being determined by a combination of Dirichlet and periodic boundary conditions. 
As Theorem~\ref{lkdvthm2} shows, the latter exhibits revivals of weak cusp type instead, and the proof of this takes into account the more complex modularity and periodicity structure of the family $\{\re^{i(2n-\frac13)\pi x}\}_{n\in\mathbb{N}}$.

\subsection*{The dislocated Dirichlet Laplacian Schr\"odinger equation}
The revival property for the boundary value problem \eqref{disPDE} is more involved than that for \eqref{pseD1}, but the result is analogous, separately in the sub-intervals $(0,b)$ and $(b,1)$. 
The next theorem gives the explicit formulation in the sub-interval $(0,b)$. We then indicate how the statement for $x\in(b,1)$ can be derived using reflections. Notably, the revival effect occurs at different times in $(0,b)$ and $(b,1)$. Moreover,  discontinuities do not propagate through the dislocation barrier. In order to facilitate the comparison with Theorem~\ref{lkdvthm2}, we have intentionally matched the notation. 

\begin{theorem} \label{Revival in (0,b)}
Let $u_0:[0,1]\longrightarrow \mathbb{C}$ be piecewise Lipschitz. Then, for all $x\in(0,b)$ and $t\in \mathbb{R}$, the solution of the boundary value problem \eqref{disPDE} admits the representation
\[
u(x,t)=U_{\mathcal{R}}(x,t)+U_{\mathcal{C}}(x,t),
\]
where $U_{\mathcal{C}}(\argdot,t)\in \ContinuousSpace(0,b)$  and
\be\label{URdis}
    U_{\mathcal{R}}(x,t)=\sum_{n=1}^\infty\widetilde{u_0}(n)\re^{-i\left(\frac \pi b(n+\frac 1 4)\right)^2t}\sin\left(\frac{\pi}{b} \Big(n+\frac14\Big) x\right)
\ee
with
\be\label{U0n}
  \widetilde{u_0}(n)=2\int_0^1
  \left[ u_0(by) \sin\Big( \pi\Big(n+\frac14\Big) y\Big)
 + u_0(b^+) \cos\Big( \pi\Big(n+\frac14\Big) y\Big)\right]\mathrm{d}y.
   \ee
If $p,q\in\mathbb{N}$ are coprime, then for all $x\in(0,b)$
$$
   U_{\mathcal{R}} \Big(x,\frac{2b^{2}}{\pi}\frac{p}{q}\Big)  =L_1(x)+L_2(x)+L_3(x)
   $$
  where
    \begin{align*}
        L_1(x)   &=\frac 1 2\sum_{k=0}^{2q-1}d_k^{p,q}\left[ i\mathrm{e}^{-i\frac{\pi }{4b}x}G_{u_0}^b\Big(-\frac{x}{b}- \frac{k}{q}\Big)-i\mathrm{e}^{i\frac{\pi }{4b}x}G_{u_0}^b\Big(\frac{x}{b}-\frac{k}{q}\Big)\right],
        \\
        L_2(x)   &=\frac 1 2\sum_{k=0}^{2q-1}d_k^{p,q}
        \left[ \mathrm{e}^{i\frac{\pi }{4b}x} [\mathcal{H}G_{u_0}^b]\Big(\frac xb -\frac kq\Big) - \mathrm{e}^{-i\frac{\pi }{4b}x} [\mathcal{H}G_{u_0}^b]\Big(-\frac xb -\frac kq\Big)  \right],
        \\
        L_3(x)   &=-\left[\int_0^2G_{u_0}^b(y)dy\right]
        \sin\Big(\frac{\pi x}{4b}\Big)
        \re^{-i\frac{\pi p}{2q}}
    \end{align*}
and
  \begin{equation}
		\label{v,g0s}
		\begin{aligned}
			&G_{u_0}^b(x) = e^{- i \frac{\pi x}{4}}
		 \begin{cases}
			u_0(b^+)+iu_0(bx), &\quad 0\leq x \leq 1, \\
			i u_0(b^+)+u_0(b(2-x)), &\quad 1\leq x \leq 2.
		\end{cases}
		\end{aligned}
	\end{equation}
Here $\mathcal{H}$ denotes the 2-periodic Hilbert transform and $d_k^{p,q}$ are complex coefficients independent of $u_0$.
\end{theorem}

The result of this theorem implies that the revival part of the solution at $x\in(0,b)$ is completely characterised by a function $U_{\mathcal{R}}$ which has an explicit expression in terms of a trigonometric series. As in Theorem~\ref{lkdvthm2}, the coefficients of this series, given in \eqref{U0n}, depend on the initial function. However,  in this case the boundary contribution is due only to the dislocation. Also, at rational times, the series becomes a finite summation of components that depend on an explicit function $G_{u_0}^b$ and on its Hilbert transform. 
In contrast to Theorem~\ref{lkdvthm2}, the term $L_3$ is not constant. However, this term is a continuous function, thus overall ensuring that \eqref{disPDE} exhibits a weak revival property.

If $u_0$ has discontinuities, the explicit formula for the revival effect in $(0,b)$ involves only the one that lie in $(0,b)$, as well as the values $u_0(0)$, $u_0(b^-)$ and $u_0(b^+)$. 
For $u_0$ compatible with the boundary and dislocation conditions, so that $u_0(0)=0$ and $u_0(b^-)=u_0(b^+)$, only the jumps of $u_0$ in $(0,b)$ will propagate as jumps and logarithmic cusps at times rationally related to $b$, and they will contribute to form the discontinuities of the solution at $(0,b)$. Indeed, the function $G^b_{u_0}$ satisfies
\begin{align*}
G^b_{u_0}(0)&=u_0(b^+)+iu_0(0) \qquad \text{and} \\
G^b_{u_0}(2)&=u_0(b^+)-iu_0(0). \end{align*}
Hence  $G^b_{u_0}$ is 2-periodic if and only if $u_0(0)=0$. If this is not the case, the function $G^b_{u_0}$ has a discontinuity at the boundary point $x=0$ that propagates into the solution. In addition, provided $u_0(b^-)=u_0(b^+)$ and therefore $u_0$ is compatible with the constraint at the dislocation, we then have
$G^b_{u_0}(1^-)=G^b_{u_0}(1^+)$. On the other hand, if $u_0$ is not continuous at $x=b$, then this additional discontinuity of $G^b_{u_0}$  also propagates into the solution.

\begin{figure}
        \centering
        \includegraphics[width=0.98\linewidth]{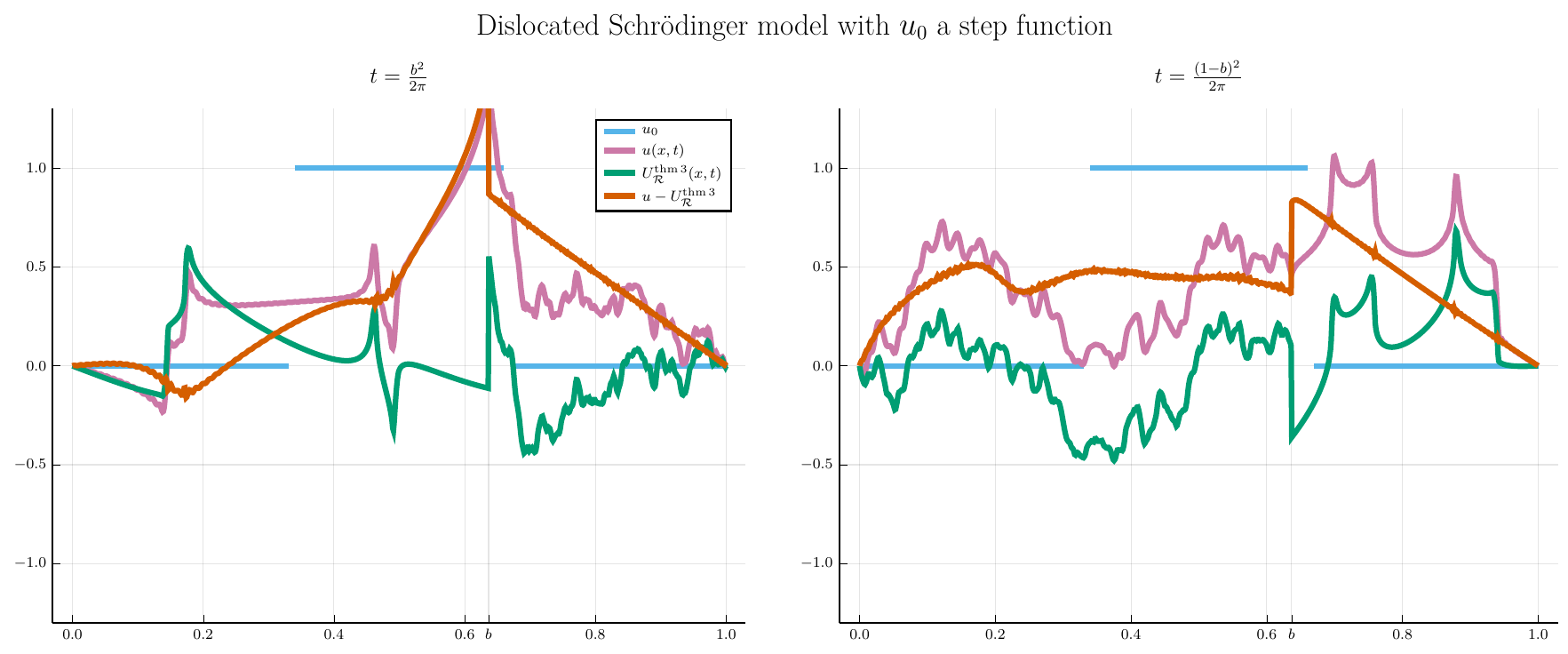}
        \caption{
 For $u_0$ a step function (blue) the graphs show the solution of \eqref{disPDE} (purple), the function $U_{\mathcal R}$ (green) and the continuous function $U_{\mathcal C}$ (brick). See Theorem~\ref{Revival in (0,b)}. The horizontal axis is $x\in[0,1]$. The time $t$ is fixed rational with respect to $(0,b)$ on the left and fixed rational with respect to $(b,1)$ on the right. These graphs were constructed from an expansion in the first 250 eigenfunctions, approximated numerically.  
        \label{fig:intro.D}}
    \end{figure}

The proof of Theorem~\ref{Revival in (0,b)} is given in Section~\ref{sec3}.   We illustrate the statement in  Figure~\ref{fig:intro.D}.  In this figure we plot the solution at times that are rational with respect to $b$ (left) and $1-b$ (right). At each of these values of the time variable, one part of the spatial interval exhibits cusp revivals, while the other part appears to be fractalised.
   
We now sketch how to formulate a statement similar to Theorem~\ref{Revival in (0,b)}, but in the interval $(1,b)$. Further details are given in Remark~\ref{reversing}. Let  $R:\Ltwozo\longrightarrow \Ltwozo$ be the reflection operator $Ru_0(x)=u_0(1-x)$. Then $R$ maps the solution of \eqref{disPDE} for $b$ into the solution of   \eqref{disPDE} for $1-b$, reflected about $\frac12$ for the time reversed to $-t$. Let $\widetilde{v_0}$ have the same expression as $\widetilde{u_0}$, with  $1-b$ replacing $b$, $Ru_0$ replacing $u_0$ and $b^-$ replacing $b^+$. Then, for $x\in(b,1)$, the solution is given by
\[
u(x,t)=V_{\mathcal{R}}(x,t)+V_{\mathcal{C}}(x,t),
\]
where $V_{\mathcal{C}}(\argdot,t)\in \operatorname{C}\big(b,1\big)$ and
\begin{equation*}
    V_{\mathcal{R}}(x,t)=\sum_{n=1}^\infty \widetilde{v_0}(n)\re^{i\left(\frac \pi b(n+\frac 1 4)\right)^2t}\sin\left(\frac{\pi}{b} \Big(n+\frac14\Big) (1-x)\right).
\end{equation*}
At the (negative) rational times,  we have
\[V_{\mathcal{R}}\Big(x,-\frac{2(1-b)^{2}}{\pi}\frac{p}{q}\Big)=M_1(x)+M_2(x)+M_3(x)\] where $M_j$ have similar expressions as $L_j$ but replacing $u_0$ by $Ru_0$, $b$ by $1-b$, and $x$ by $1-x$. In this case, only the discontinuities of $u_0$ in the interval $(b,1)$ propagate, as jumps and logarithmic cusps.

\medskip

The discussion and proof of the main results pertaining the boundary value problems~\eqref{pseD1} and \eqref{disPDE} is given in sections~\ref{sec4} and \ref{sec3}. These two sections comprise the technical part of this work. In Section~\ref{rev} we describe a general framework fitting the new manifestation of the revival property described in theorems \ref{lkdvthm2} and \ref{Revival in (0,b)}. Appendix~\ref{apa} includes details of the main definitions and notation that we use throughout the work. The final Appendix~\ref{apb} is devoted to the operator-theoretical properties of the spatial operators.


\section{Proof of Theorem~\ref{lkdvthm2}} \label{sec4}

The boundary value problem \eqref{pseD1} is the time-evolution equation for the third-order differential operator $iA$, where $A:\operatorname{Dom}(A)\longrightarrow L^2(0,1)$ is given by
$$
A\phi(x)=i\phi_{xxx}(x)
$$
and
\be\label{airydom1}
\operatorname{Dom}(A)=\big\{\phi\in H^3(0,1):\;\phi(0)=0,\;\phi(1)=0,\; \phi_x(1)= \phi_x(0)\big\}.
\ee
 The linear operator $A$ is self-adjoint and it has a compact resolvent, see Appendix~\ref{apb}. Hence, for all $u_0\in\Ltwozo$, we have that \eqref{pseD1} admits a unique solution which has a series expansions in terms of the orthonormal basis of eigenfunctions of $A$. Our strategy for the proof of Theorem~\ref{lkdvthm2} is to determine the leading order terms of this series and then show that these terms give rise to the cusp revival effect.

\medskip

We first describe the eigenvalues and eigenfunctions of the operator $A$. The eigenvalues are given by $\lambda=k^3$, where $k$ are the nonzero roots\footnote{Although $0$ is a root of $\Delta$, notice that no linear combination of polynomials of order 2, other than the zero polynomial, satisfy all boundary conditions simultaneously.} of the characteristic determinant
$\Delta(k)=0$, given by
\begin{align}
\Delta(k)&=\re^{ik}+\re^{- ik}+\alpha(\re^{ i\alpha k}+\re^{-i\alpha k})+\alpha^2(\re^{i\alpha^2 k}+\re^{-i\alpha^2 k}),\nonumber \\
\alpha&=- \frac 1 2 +i\frac{\sqrt 3}2=\re^{\frac{2\pi i}3}. \label{DeltaAlpha}
\end{align}
Below we often use the fact that $\alpha$ is a cubic roots of unit, so that 
$\overline{\alpha}=\alpha^2=\alpha^{-1}$.

 The spectrum of $A$ is a countable set of real eigenvalues, each of finite multiplicity, accumulating at $\pm \infty$. The zeros of $\Delta$ lie on the lines $\mathbb{R}$, $\alpha \mathbb{R}$ and $\alpha^2 \mathbb{R}$, see~\cite{pelloni2005spectral}.
Moreover $\Delta(\pm\alpha k) = \alpha^2\Delta(k)$ for all $k\in\mathbb{C}$.
Therefore, the nonzero roots of $\Delta$ are of the form
$\pm\alpha^jk_n$ for $j\in\{0,1,2\}$ where $\{k_n\}_{n\in\N}$ is the increasing sequence of positive real roots of $\Delta$.
Thus
\[
   \operatorname{Spec}(A)=\big\{\lambda_n\big\}_{n\in\mathbb{Z}\setminus\{0\}},
\]
where $\lambda_n= k_n^3$ for $n\in\N$ and $\lambda_{-n}=-\lambda_n$. As we shall see next, the eigenvalues of $A$ are all simple. The zeros $k_n$ of $\Delta$ cannot be computed in closed form, however we can characterise their asymptotic behaviour.

\begin{prop} \label{prop:Airy.eigenvalues}
The positive real roots $k_n$ of the equation $\Delta(k)=0$ are all simple and have the following asymptotic behaviour,
\be\label{evas}
k_n=\left(2n-\frac 1 3 \right)\pi+\bigoh{\re^{-\sqrt{3} n \pi}} \text{ as } \ n\to\infty.
\ee
\end{prop}
\begin{proof}
Rewrite the determinant as
$$
\Delta(k)=2\left[\cos(k)-\cos\left(\frac k 2\right)\cosh\left(\frac{\sqrt 3 k}2\right)-
\sqrt 3\sin\left(\frac k 2\right)\sinh\left(\frac{\sqrt 3 k}2\right)\right].
$$
Then, $\Delta(k) \in\R$ for all $k\in\R$. By rearranging this, assuming that $k\not=0$, we find that any real zero of $\Delta$  must satisfy
\be\label{RHS}
\tan\left(\frac k 2\right)=\frac{\cos\left(\frac k 2\right)-\cosh\left( \frac {\sqrt 3} 2 k\right)}{\sin\left(\frac k 2\right)+\sqrt 3\sinh\left( \frac {\sqrt 3} 2 k\right)}.
\ee
As $k\to\infty$, the limit of the right of equation~\eqref{RHS} is $-\frac{1}{\sqrt{3}}$.
This, and a routine computation, give the asymptotic expression~\eqref{evas}.

Now, we show that all the positive roots $k_n$ are simple zeros. Since the right side of~\eqref{RHS} is negative everywhere, all $k_n$ must lie in intervals $((2n-1)\pi,2n\pi)$; in particular, there is no positive root $k\leq\pi$.
The left side of equation~\eqref{RHS} has gradient at least $\frac12$ everywhere.
To complete the proof that all of the $k_n$ are simple,
it suffices to show that the right side of~\eqref{RHS} has gradient less than $\frac12$ everywhere in $(\pi,\infty)$. This follows from a direct calculation, which we omit.
\end{proof}

By combining this proposition with~\cite[Theorem~4.1]{locker}, it follows that all the eigenvalues of $A$ are simple.

\medskip

Let us now consider the eigenfunctions of $A$. The structure of these eigenfunctions was originally examined in \cite{pelloni2005spectral, pelloni2013spectral}. We fix the notation as follows. Let  
\begin{equation} \label{eqn:Airy.ef.npositive}
\phi_{n}(x)=\sum_{r=0}^2\left[ \re^{i\alpha^{r+2}k_n} - \re^{i\alpha^{r+1}k_n}\right]\re^{i\alpha^{r}k_nx}.
\end{equation}
Then $\phi_n(x)$ is an eigenfunction, corresponding to the eigenvalue $k_n^3$. 
For each $n\in\N$, there is an eigenfunction $\phi_{-n}$ corresponding to $-k_n^3$, which is given by replacing $k_n$  with $-k_n$ on the right side of~\eqref{eqn:Airy.ef.npositive}.
Note that
    \begin{equation} \label{normsq}
        \lVert\phi_n\rVert^2 = \re^{\sqrt{3}k_n}\left(1+\frac 1 {\sqrt 3k_n}\right)+\bigoh{k_n^{-1}\re^{\sqrt3k_n/2}}\;\;
   \text{ as } \ n\to\infty.
    \end{equation}
The family
\[
     \left\{\|\phi_n\|^{-1}\phi_n\right\}_{n\in\mathbb{Z}\setminus\{0\}}\subset \mathrm{L}^2(0,1)
\]  
is an orthonormal basis of eigenfunctions.

In our arguments below, it is convenient to combine the contributions of the two eigenfunctions $\phi_{\pm n}$ into a single term. Since
$$
\phi_{-n}(x)=-\overline{\phi_n(x)},
$$
then for any $u_0\in L^2(0,1)$, the solution to the boundary value problem \eqref{pseD1} is given by
\begin{align}\notag
u(x,t) &=
\sum_{n=1}^\infty \left[\frac{\langle u_0,\phi_n\rangle}{\|\phi_n\|^2}\re^{ik_n^3t}\phi_n(x)+\frac{\langle u_0,\phi_{-n}\rangle}{\|\phi_{-n}\|^2}\re^{-ik_n^3t}\phi_{-n}(x)\right] \\
\label{urep2}
&= \sum_{n=1}^\infty 2\Re \left[\frac{\langle u_0,\phi_n\rangle}{\|\phi_n\|^2}\re^{ik_n^3t}\phi_n(x)\right].
\end{align}   

\medskip
    
We now establish the asymptotic behaviour of the terms
\[
   \frac{\langle u_0,\phi_n\rangle}{\|\phi_n\|^2} \phi_n(x),
\] for $u_0$ satisfying the hypothesis of Theorem~\ref{lkdvthm2}.

\begin{prop} \label{prop:Airy.simplifiedeigenfunctions}
 Let $u_0:[0,1]\longrightarrow \mathbb{R}$ be piecewise Lipschitz. Then,
    \begin{equation*}
        \frac{\left\langle u_0,\phi_n\right\rangle}{\left\lVert\phi_n\right\rVert^2} \phi_n(x)
        =
        \re^{ik_nx} \left[ \left\langle u_0 , \re^{i k_n (\argdot)} \right\rangle + \frac{i\alpha}{k_n} \left(u_0(1)+u_0(0)\right) \right] + \bigoh{n^{-\frac32}}
    \end{equation*}
as $n\to\infty$,  uniformly on compact subsets for all $x\in (0,1)$. 
\end{prop}
\begin{proof}
Firstly note that $k_n^{-1}\sim \frac{1}{2\pi} n^{-1}$, 
\[
         \re^{\pm i\alpha k_n} \sim \re^{\mp \sqrt3 \pi n} \quad \text{and} \quad  \re^{\pm i\alpha^2k_n} \sim \re^{\pm \sqrt3 \pi n},
\]
as $n\to \infty$. Also $k_n\in \mathbb{R}$, so that $|\re^{ik_nx}|=1$. Below we will use all this repeatedly.

Let
\[
     \frac{\left\langle u_0 , \phi_n \right\rangle}{\left\lVert\phi_n\right\rVert^2} \phi_n(x)
        =
        \left\langle u_0 , \re^{i k_n (\argdot)} \right\rangle \re^{ik_nx}+A_n(x)+B_n(x)+C_n(x),
\]
where
\begin{align*}
    A_n(x)&=\frac{\displaystyle\sum_{r=1}^2\left[ \re^{-i\alpha^{-(r+2)}k_n} - \re^{-i\alpha^{-(r+1)}k_n} \right] \left\langle u_0 , \re^{i\alpha^rk_n(\cdot)} \right\rangle}{\left\lVert\phi_n\right\rVert^2}\left( \re^{i\alpha^2k_n} - \re^{i\alpha k_n} \right) \re^{ik_nx}, \\
    B_n(x)&=\left[ \frac{\left\lvert \re^{i\alpha^2k_n} - \re^{i\alpha k_n} \right\rvert^2}{\left\lVert\phi_n\right\rVert^2} - 1 \right] \left\langle u_0 , \re^{i k_n (\argdot)} \right\rangle \re^{ik_nx} \qquad \text{and}\\
    C_n(x)&=\frac{\displaystyle\sum_{r=0}^2\left[ \re^{-i\alpha^{-(r+2)}k_n} - \re^{-i\alpha^{-(r+1)}k_n} \right] \left\langle u_0 , \re^{i\alpha^rk_n(\cdot)} \right\rangle}{\left\lVert\phi_n\right\rVert^2} \sum_{j=1}^2\left( \re^{i\alpha^{j+2}k_n} - \re^{i\alpha^{j+1}k_n} \right) \re^{ik_n\alpha^jx}.
\end{align*}
The proof of the proposition reduces to showing that, as $n\to\infty$ and for all $x\in(0,1)$,
\begin{equation} \label{asympA}
       A_n(x)= \frac{i\alpha}{k_n} \left(u_0(1)+u_0(0)\right) \re^{ik_nx} + \bigoh{n^{-\frac32}},
\end{equation}
\[
B_n(x)=\bigoh{n^{-2}}\quad\text{and}\quad C_n(x)=\bigoh{\re^{-\sqrt3 x k_n/2}}.\]  
In the argument below it will be clear that the estimates for $A_n(x)$ and $B_n(x)$ hold uniformly on compact sets, but we will give full details in the case of $C_n(x)$.   We split the rest of the proof into 4 steps. 

\underline{Step~1}: to find preliminary estimates on the inner products, we first decompose $u_0$ into an absolutely continuous function $u^{\mathrm{ac}}_0$ and a piecewise constant function $u_0^j$ with jumps at finitely many point $x_j\in(0,1)$, see Remark~\ref{crucial}. 
Then we have
\begin{equation} \label{L-Stil}
\begin{aligned}
      \langle u_0,\mathrm{e}^{i\alpha^r k_n(\cdot)}\rangle&=
      \int_0^1 u_0(y)\mathrm{e}^{-i\overline{\alpha^r}k_n y}\mathrm{d}y \\
      &=\frac{i\alpha^r}{k_n}\left[u_0(1)\mathrm{e}^{-i\alpha^{-r}k_n}-u_0(0)-\int_0^1\mathrm{e}^{-i\alpha^{-r}k_ny}\mathrm{d}u_0(y)\right],     
\end{aligned}\end{equation}
where the integral on the right hand side can be decomposed as
\begin{equation} \label{L-Still2}
    \int_0^1\mathrm{e}^{-i\alpha^{-r}k_ny}\mathrm{d}u_0(y)=\int_0^1\mathrm{e}^{-i\alpha^{-r}k_ny}\mathrm{d}u^{\mathrm{ac}}_0(y)+\sum_{j=1}^k h_j \mathrm{e}^{-i\alpha^{-r}k_n x_j},
\end{equation}
with $h_j$ the jumps of $u_0$ at the points $x_j$. 

Note that, since $(u^{\mathrm{ac}}_0)'\in \mathrm{L}^\infty[0,1]$, we have
\[
   \left|\int_0^1\mathrm{e}^{-i\alpha^{-r}k_ny}\mathrm{d}u^{\mathrm{ac}}_0(y)\right|\leq  \|(u^{\mathrm{ac}}_0)'\|_{\infty}\|\mathrm{e}^{-i\alpha^{-r}k_n(\cdot)}\|_{\infty}.
\]
Hence, we have the following  estimate, used in step 3 and step 4 below. 
   \begin{equation} \label{eqn:Airy.prop2proof.extraterms.jnonzero}
        \left\lvert \left\langle u_0 , \re^{i\alpha^rk_n(\argdot)} \right\rangle \right\rvert
        =
        \begin{cases}
            \bigoh{n^{-1}} & \mbox{if } r = 0,\,1 \\
            \bigoh{\frac{\re^{\sqrt3 \pi n}}{n}} & \mbox{if } r = 2.
        \end{cases}
    \end{equation}

\underline{Step~2}: consider \eqref{asympA} for the term $A_n(x)$.  Firstly observe the reduction
 \begin{multline*} A_n(x)= 
        \frac{\re^{ik_nx}}{\left\lVert\phi_n\right\rVert^2}
        \left[
            \left\langle u_0 , \re^{i\alpha k_n (\argdot)} \right\rangle
            \left( - \re^{\sqrt3k_n} + \bigoh{\re^{\sqrt3k_n/2}} \right)
            \right. \\ \left.
            +
            \left\langle u_0 , \re^{i\alpha^2 k_n (\argdot)} \right\rangle
            \left( - \re^{k_n(\sqrt3-3i)/2} + \bigoh{1} \right)
        \right]
        \\
        = - \left\langle u_0 , \re^{i\alpha k_n (\argdot)} + \re^{ik_n} \re^{-i\alpha^2k_n(1-(\argdot))} \right\rangle \re^{ik_nx}
        + \bigoh{\frac1{n^2}}.
    \end{multline*}
We can derive a more detailed estimate than \eqref{eqn:Airy.prop2proof.extraterms.jnonzero} for the inner product appearing in the last line above. Our starting point are the expressions \eqref{L-Stil} and \eqref{L-Still2} for $r=1$. Since $(u^{\mathrm{ac}}_0)'\in \mathrm{L}^2(0,1)$, the Riemann-Lebesgue lemma implies that
\[
   \left| \int_0^1\mathrm{e}^{-i\alpha^{-1} k_ny}\mathrm{d}u^{\mathrm{ac}}_0(y)\right|=\left| \int_0^1\mathrm{e}^{i\frac{ k_n}{2}y}\mathrm{e}^{-\frac{ \sqrt{3}k_n}{2}y} (u^{\mathrm{ac}}_0)'(y) \mathrm{d}y\right|=\mathrm{o}\left(n^{-\frac{1}{2}}\right),
\]
thus
\[
      \left\langle u_0 , \re^{i\alpha k_n(\argdot)} \right\rangle =-\frac{i\alpha}{k_n}u_0(0)+\bigoh{n^{-\frac{3}{2}}}.
\]
By an analogous calculation, using that, by proposition~\ref{prop:Airy.eigenvalues}, both
    \(\alpha \re^{-i\alpha^2k_n}\)
    and
    \(\alpha^2\re^{-ik_n(1-\alpha)}\)
    are exponentially decaying, and
    \(\alpha^2 \re^{-ik_n} = -\alpha + \bigohsymbolonly(\re^{-\sqrt3n\pi})\), we have that
\[
      \left\langle u_0 , \re^{ik_n}\re^{-i\alpha^2 k_n(1-(\argdot))} \right\rangle =-\frac{i\alpha}{k_n}u_0(1)+\bigoh{n^{-\frac{3}{2}}}.
\]   
Therefore, \eqref{asympA} holds true.

\underline{Step 3}: now consider the estimate for $B_n(x)$. Since
  \[
    \frac{\left\lvert \re^{i\alpha^2k_n} - \re^{i\alpha k_n} \right\rvert^2}{\left\lVert\phi_n\right\rVert^2} - 1 =    \frac{1}{\left\lVert\phi_n\right\rVert^2} \left\lvert \re^{-ik_n} \right\rvert \left\lvert \re^{\sqrt3k_n} + \re^{-\sqrt3k_n} - 2 \right\rvert - 1 = \bigoh{\frac1{k_n}},
    \] 
    then, on account of \eqref{eqn:Airy.prop2proof.extraterms.jnonzero}, we have $B_n(x)=\bigoh{n^{-2}}$.
    
\underline{Step~4}: finally consider the estimate for $C_n(x)$. 
Note that,
    \[
        \sum_{j=1}^2\left( \re^{i\alpha^{j+2}k_n} - \re^{i\alpha^{j+1}k_n} \right) \re^{ik_n\alpha^jx}
        =
        - \re^{i\alpha^2 k_n} \re^{ik_n\alpha x} - \re^{i k_n} \re^{ik_n\alpha^2 x}
        +
        \bigoh{1},
    \]
    uniformly in $x\in[0,1]$.
    Fixing a small $\epsilon>0$, we obtain the tighter estimate
    \[
        \sum_{j=1}^2\left( \re^{i\alpha^{j+2}k_n} - \re^{i\alpha^{j+1}k_n} \right) \re^{ik_n\alpha^jx}
        =
        \bigoh{\re^{\sqrt3(1-\epsilon)k_n/2}},
    \]
    uniformly in $x\in[\epsilon,1-\epsilon]$.
    The asymptotic formula~\eqref{eqn:Airy.prop2proof.extraterms.jnonzero} implies that
    \[
        \sum_{r=0}^2\left[ \re^{-i\alpha^{-(r+2)}k_n} - \re^{-i\alpha^{-(r+1)}k_n} \right] \left\langle u_0 , \re^{i\alpha^rk_n} \right\rangle
        =
        \bigoh{\re^{\sqrt3k_n/2}}.
    \]
    Combining the last two asymptotic bounds with estimate~\eqref{normsq} in the formula for $C_n(x)$, yields \[C_n(x)=\bigohsymbolonly(\re^{-\sqrt3\epsilon k_n/2})\] uniformly for $x\in[\epsilon,1-\epsilon]$. This completes the proof of the proposition.    
\end{proof}

\medskip

In the proof of the first part of Theorem~\ref{lkdvthm2}, we will construct the continuous perturbation $U_{\mathcal C}$ as
$$
    U_{\mathcal C}(x,t)=w_1(x,t)+w_2(x,t),
$$
where $w_1,w_2$ are continuous functions of $x$ on $(0,1)$. For that purpose we will argue as follows. 
Starting with the solution representation formula~\eqref{urep2}, we replace the normalised eigenfunctions with their approximations $\re^{ik_nx}$. If $u_0$ does not vanish at the boundary, then, as stated in Proposition~\ref{prop:Airy.simplifiedeigenfunctions}, we also include the boundary values. The function  $w_1$ will be defined as  the difference between the original solution representation~\eqref{urep2} and this approximate representation. In the expression for $u-w_1$ we replace $k_n$ with their leading order terms $(2n-\frac13)\pi$, thus obtaining the formula for $U_{\mathcal{R}}$.
    The function  $w_2$ will then be defined as the difference between $u-w_1$ and $U_{\mathcal{R}}$.

\medskip

Before proceeding to the proof of Theorem~\ref{lkdvthm2}, we establish a technical lemma. Let $s\in\mathbb{R}$ be a fixed parameter. We define the \emph{periodic translation} operator $T_s:\Ltwozo\longrightarrow \Ltwozo$ as
\be
\label{transl}
\quad T_s v(x) = v^{\dagger}(x-s),
\ee
where  $v^{\dagger}$ is the $1$-periodic extension of $v$ from $[0,1)$ to $\R$.
For $n\in\Z$, we will write $e_n(x) = \re^{2\pi i nx}$.

\begin{lemma} \label{lem:Airy.MainRevivalLemma}
    For  $p,q\in\N$ coprime, let
    \be\label{dkpq}
        d_k^{p,q} = \frac1q \sum_{m=0}^{q-1} \re^{2\pi i\left( mk + \left[4m^3-2m^2\right]p \right)/q}.
    \ee
    Then, for all $v\in \Ltwozo$ we have
    \[
        \sum_{n=1}^\infty
        \left\langle v , e_n \right\rangle
        \re^{i(2n)^3\pi \frac pq} \re^{-i (2n)^2\pi\frac pq}
        e_n(x)
        =
        \sum_{k=0}^{q-1} d_k^{p,q} \, T_{\frac kq} \sum_{n=1}^\infty
        \left\langle v , e_n \right\rangle e_n(x).
    \]
\end{lemma}

\begin{proof}
    Fix $p,q\in\N$ coprime.
    If $m,n\in\N$ with $m \equiv_q n$ then, $m^r \equiv_q n^r$ for all $r\in \mathbb{N}$. Thus, $4m^3-2m^2 \equiv_q 4n^3-2n^2$. Hence,
    \[
        \re^{2\pi i (4m^3-2m^2) \frac pq} = \re^{2\pi i (4n^3-2n^2) \frac pq}.
    \]
    Because it represents a sum of roots of unity,
    \be\label{unitroots}
        \sum_{k=0}^{q-1} \re^{2\pi i(m-n)\frac kq} =
        \begin{cases}
            q & \mbox{if } m \equiv_q n \\
            0 & \mbox{otherwise.}
        \end{cases}
    \ee
    Summing over $m\in\{0,1,\ldots,q-1\}$ and dividing by $q$, gives
    \begin{equation} \label{eqn:RevivalMechanismGeneral}
        \sum_{k=0}^{q-1} d_k^{p,q} \re^{-2\pi i \frac{nk}q}
        = \frac1q \sum_{m=0}^{q-1}\left[\sum_{k=0}^{q-1}  \re^{2\pi i \frac{mk}q}\re^{-2\pi i \frac{nk}q} \right] \re^{2\pi i (4m^3-2m^2) \frac pq}
        = \re^{2\pi i (4n^3-2n^2) \frac pq}.
    \end{equation}
Note that, this representation allows us to replace higher powers of $n$ on the right, with only a first power of $n$ on the left, at the cost of a finite sum.

Substituting equation~\eqref{eqn:RevivalMechanismGeneral}
into the expression
\[
        \sum_{n=1}^\infty
        \left\langle v , e_n \right\rangle
        \re^{i(2n)^3\pi \frac pq} \re^{-i (2n)^2\pi\frac pq}
        e_n\left(x\right),
\]
yields
    \begin{align*}
        \sum_{k=0}^{q-1} d_k^{p,q} \sum_{n=1}^\infty
        \left\langle v , e_n \right\rangle \re^{-2\pi i \frac{nk}q} e_n(x)
        &=
        \sum_{k=0}^{q-1} d_k^{p,q} \, T_{\frac kq} \sum_{n=1}^\infty
        \left\langle v , e_n \right\rangle e_n(x),
    \end{align*}
    as required.
\end{proof}

\medskip

\begin{proof}[Proof of Theorem~\ref{lkdvthm2}]
We split the proof into two steps.

\underline{Step 1}: proof of the first part. Let
    \begin{align*}
        &w_1(x,t) = u(x,t) - \sum_{n=1}^\infty 2\Re \left[  \left( \langle u_0,\re^{ik_n(\argdot)}\rangle + \frac{i\alpha}{k_n}\left[u_0(1)+u_0(0)\right] \right) \re^{ik_nx}\re^{ik_n^3t} \right] \\
        &= \sum_{n=1}^\infty 2\Re \left[   \frac{\left\langle u_0 , \phi_n \right\rangle}{\left\lVert\phi_n\right\rVert^2}\phi_n(x) - \left( \langle u_0,\re^{ik_n(\argdot)}\rangle + \frac{i\alpha}{k_n}\left[u_0(1)+u_0(0)\right] \right) \re^{ik_nx}\re^{ik_n^3t}  \right]
    \end{align*}
 According to Proposition~\ref{prop:Airy.simplifiedeigenfunctions}, the $n$-th term in this series is $\bigohsymbolonly(n^{-\frac32})$ uniformly in $x$ on compact sub-sets of $(0,1)$.
Since the series is uniformly absolutely convergent in $x$ on compact sub-sets of $(0,1)$ and each term in it is continuous in $x$, then $w_1$ is continuous in $x$.
    
    Let
    \begin{equation} \label{lambdan}
        \kappa_n = \left(2n-\frac13\right)\pi.
    \end{equation}
Let 
\[
        g_n(x) = \left( \re^{i(k_n-\kappa_n)x}-1 \right) \qquad \text{and} \qquad h_n(t) = \left( \re^{i(k_n^3-\kappa_n^3)t}-1 \right) \re^{i\kappa_n^3t}.
    \]
   Since
    \[
        k_n-\kappa_n = \bigoh{\re^{-\sqrt3n\pi}} \qquad \text{and} \qquad k_n^3-\kappa_n^3 = \bigoh{n^2\re^{-\sqrt3n\pi}},
    \]
 then, $g_n(x) = \bigohsymbolonly(\re^{-\sqrt3n\pi})$ uniformly for $x\in [0,1]$ and $h_n(t)=\bigohsymbolonly(n^2\re^{-\sqrt{3}n\pi})$ uniformly for $t\in[0,\tau]$, where  $\tau$ is fix but arbitrarily large.

    Now, let
    \begin{align*}
        f_n &= \langle u_0 , \re^{ik_n(\argdot)} \rangle - \langle u_0 , \re^{i\kappa_n(\argdot)} \rangle + i\alpha\left(\frac1{k_n}-\frac1{\kappa_n}\right) \left[ u_0(1) + u_0(0) \right]
        \\
        &= \int_0^1 u_0(x) \re^{-i\kappa_n x}\left( \re^{-i(k_n-\kappa_n)} - 1 \right) \mathrm{d}x + i\alpha\left(\frac{\kappa_n-k_n}{k_n\kappa_n}\right) \left[ u_0(1) + u_0(0) \right].
    \end{align*}
By virtue of Proposition~\ref{prop:Airy.eigenvalues}  and the above estimate on $\kappa_n-k_n$, then the second term in this expression is \[ i\alpha\left(\frac{\kappa_n-k_n}{k_n\kappa_n}\right) \left[ u_0(1) + u_0(0) \right]=\bigohsymbolonly(n^{-2}\re^{-\sqrt3n\pi})\]
as $n\to \infty$.
 Moreover, the modulus of the integral in the first term, is  bounded above by
    \[
        \left\lVert u_0 \right\rVert_\infty \int_0^1 \left\lvert \re^{-i(k_n-\kappa_n)x} - 1 \right\rvert \mathrm{d}x
        \leq \left\lVert u_0 \right\rVert_\infty \int_0^1 \left\lvert (k_n-\kappa_n)x \right\rvert \mathrm{d}x
        = \bigoh{\re^{-\sqrt3n\pi}}.
    \]
    
Define $w_2=u-w_1-U_{\mathcal{R}}$.
    Note that
    \[
        \int_0^1 \re^{-i\left(2n-\frac13\right)\pi y}\mathrm{d} y =\frac{i}{\left(2n-\frac13\right)\pi}(\re^{i\pi/3}-1)= \frac{i\alpha}{\kappa_n},
    \]
    so the coefficients $\widetilde{u_0}(n)$ in the expression \eqref{eqn:defn.u0hat} are given by
    \[
         \widetilde{u_0}(n) =\left\langle u_0+u_0(1)+u_0(0),\re^{i\kappa_n\pi(\argdot)}\right\rangle=
        \left\langle u_0 , \re^{i\kappa_n\pi(\argdot)} \right\rangle
        + \frac{i\alpha}{\kappa_n} \left[ u_0(1) + u_0(0) \right].
    \]
Substituting this in the representation of $U_\mathcal{R}$ and reducing the three term sum inside each term of the series, gives 
    \begin{multline*}
        w_2(x,t) = \sum_{n=1}^\infty 2 \Re\Bigg[
            h_n(t) \left( \langle u_0,\re^{ik_n(\argdot)}\rangle + \frac{i\alpha}{k_n}\left[u_0(1)+u_0(0)\right] \right) \re^{ik_nx} \\
            + \re^{i\kappa_n^3t} f_n \re^{ik_nx} 
            + \re^{i\kappa_n^3t} \left( \langle u_0,\re^{i\kappa_n(\argdot)}\rangle + \frac{i\alpha}{\kappa_n}\left[u_0(1)+u_0(0)\right] \right) g_n(x)
        \Bigg].
    \end{multline*}
    
    Now, since $k_n,\,\kappa_n\in \mathbb{R}$, it follows that the terms multiplying each of $h_n(t)$, $f_n$ and $g_n(x)$ are bounded uniformly in $n\in\N$, for $x\in[0,1]$ and $t\in[0,\tau]$.
Hence, by the asymptotic bounds obtained above for $h_n(t)$, $f_n$ and $g_n(x)$, all the terms in the series for $w_2$ are $\bigohsymbolonly(n^2\re^{-\sqrt{3}n\pi})$ as $n\to \infty$, uniformly in $x\in[0,1]$ and $t\in[0,\tau]$.
 Thus, by Weierstrass' M-test, $w_2$ is continuous.

 Since $w_1+w_2 = u - U_{\mathcal{R}}$ this completes the proof of the first part of the theorem.

\medskip

\underline{Step 2}: we consider the second part of the  theorem.  Let    \[
        U_{\mathcal R}(x,t)
        =
        \sum_{n=1}^\infty 2\Re\left[\left\langle u_0(\argdot)+u_0(1)+u_0(0),\re^{i(2n-\frac 13)\pi \argdot}\right\rangle\re^{i(2n-\frac 1 3)^3\pi^3 t}\re^{i(2n-\frac 13)\pi x}\right].
    \]
 Since
    \[
        \re^{i(2n-\frac 13)\pi x}
        = \re^{2\pi i nx}\re^{-i\frac{\pi x}3}
        = e_n(x)\re^{-i\frac{\pi x}3},
    \]
    expanding the cubic power in the time exponential $\re^{i(2n-\frac 1 3)^3\pi^3 t}$, yields
       \[
        U_{\mathcal R}(x,t)
        =
        2\Re\left[
            \re^{-i\frac{\pi }3 x}\re^{-i(\frac\pi 3)^3t}\sum_{n=1}^\infty \left\langle G_{u_0}, e_n \right\rangle  \re^{i(2\pi n)^3t}\re^{2\pi i n\frac {\pi^2}3 t}\re^{-i (2\pi n)^2\pi t} e_n(x)
        \right],
    \]
    where $G_{u_0}$ is given by \eqref{v0NEW}.
    Note that the series in the expression above
    is analogous to the one for the quasi-periodic Airy problem, studied in \cite{BFP}, except that the summation is only over $n\geq 1$.
   
Now,
    \[
        U_{\mathcal R}\left(x,\frac p{q\pi^2}\right)
        =
        2\Re\left[
            \re^{-i\frac{\pi }3 x} \re^{-i \frac\pi {27}\frac p q}
            T_{\frac{-p}{3q}}
            \sum_{n=1}^\infty
            \left\langle G_{u_0} , e_n \right\rangle
            \re^{i(2n)^3\pi \frac pq} \re^{-i (2n)^2\pi\frac pq}
            e_n(x)
        \right].
    \]
 According to Lemma~\ref{lem:Airy.MainRevivalLemma}, this can be expressed also as
    \[
        U_{\mathcal R}\left(x,\frac p{q\pi^2}\right)
        =
        2\Re\left[
            \re^{-i\frac{\pi }3 x} \re^{-i \frac\pi {27}\frac p q}
            \sum_{k=0}^{q-1} d_k^{p,q} \, T_{\frac kq - \frac p{3q}}
            \sum_{n=1}^\infty
            \left\langle G_{u_0} , e_n \right\rangle
            e_n(x).
        \right]
    \]
By virtue of a canonical identity for the Hilbert transform, see \eqref{HT identities}, we have
    \[
        2\sum_{n=1}^\infty \left\langle G_{u_0},e_n\right\rangle e_n(x)
        =
        G_{u_0}(x) + i[{\mathcal H}G_{u_0}](x) - \int_0^1 G_{u_0}(y) \mathrm{d}y.
    \]
Then, since \[\sum_{k=0}^{q-1}d_k^{p,q}=1,\] the conclusion of the theorem follows.
\end{proof}


\section{Proof of Theorem~\ref{Revival in (0,b)}} \label{sec3}

Let $b\in(0,1)$ be a fixed parameter. The boundary value problem \eqref{disPDE} is the time-evolution equation associated with the linear operator $iD$, where
\begin{equation} \label{expressiond}
   D= \partial_x\big(\text{sgn}(\argdot-b) \partial_x\big):\mathrm{Dom}(D)\longrightarrow \Ltwozo.
\end{equation}
We define the domain $\operatorname{Dom}(D)$ of this quasi-differential expression as the set of all
$u:[0,1]\longrightarrow \mathbb{C}$ such that
\begin{itemize}
    \item
    $u(x)$ is absolutely continuous in $[0,1]$
    \item
    $\operatorname{sgn}(x-b)u'(x)$ is absolutely continuous in $[0,1]$
    \item and
    $u(0)=u(1)=0$.
\end{itemize}
Note that any $u\in\operatorname{Dom}(D)$ is continuous and its derivative  has a discontinuity at $x=b$ with
\begin{equation} \label{inter}
 u_x(b^+) =-u_x(b^-).
\end{equation}
This interface condition is crucial in what follows.

We associate with $D$ the boundary value problem for the linear Schr\"odinger equation given by
\begin{equation*}
\label{Dislocation model PDE}
    \partial_{t}u = iDu, \quad u|_{t=0}=u_{0}\in \Ltwozo.
\end{equation*}
See \cite{JussiDavid2018}.
Written out explicitly, this boundary value problem is given by
\eqref{disPDE}. We show that the effect of the dislocation at $x=b$ results in a  solution structure that shares important similarities with the one described for \eqref{pseD1}. In particular, we prove that, if $u_0$ has  jump discontinuities, the solution displays both discontinuities and logarithmic cusps,  characterised by the periodic Hilbert transform of functions associated to $u_0$. However, the times at which these occur and their position are different in the sub-intervals $(0,b)$ and $(b,1)$, unless $b=\frac12$.

\medskip

In Appendix~\ref{apb}, we show that $D$ is a self-adjoint operator with compact resolvent. Therefore, $D$ has purely discrete spectrum accumulating at $\pm\infty$. We now find the eigenvalues and eigenfunctions of $D$, and show that the eigenvalues have a quadratic leading order while the eigenfunctions are of an exponential type with fast decaying  reminders.  This is the first step in the proof of Theorem~\ref{Revival in (0,b)}.

The eigenvalue problem associated to $D$ is
\begin{equation} \label{dislo2}
\begin{cases} \phi_{xx}(x)+\la \phi(x)=0 & 0\leq x < b \\ \phi_{xx}(x)-\la \phi(x)=0 & b<x\leq 1 \\
 \phi(0)=\phi(1)=0, \end{cases}
\end{equation}
subject to continuity at $b$ and the interface condition \eqref{inter}, which we re-write as
\begin{equation} \label{dislo1}
\phi(b^-)=\phi(b^+) \qquad \text{and} \qquad \phi_x(b^-)=-\phi_x(b^+).
\end{equation}
Let $\alpha,\beta\in\mathbb{C}$ be constants. For $\la=0$, the general solution to \eqref{dislo2} is
\begin{equation}\label{p-wlinear} \phi(x)=\begin{cases}  \alpha x & 0\leq x < b \\ \beta (1-x) & b < x \leq 1. \end{cases}\end{equation}
If $\la=k^2\not=0$ for $k\in\mathbb{C}$, the the solutions are
\begin{equation} \label{p-wtrig} \phi(x)=\begin{cases}\alpha \sin(k x) & 0\leq x<b \\ \beta \sinh(k(1-x)) & b<x \leq 1. \end{cases}\end{equation}
Here $k\in\C$ is an eigenvalue
whenever $\alpha,\,\beta$ and $k$ are such that $\phi\in\mathrm{Dom}(D)\setminus\{0\}$, in particular whenever $\phi(x)$ also satisfies \eqref{dislo1}.

From the expressions \eqref{p-wlinear} and \eqref{p-wtrig}, for the general solutions to the eigenvalue problem associated to $D$, it follows that all the eigenvalues of $D$ must be simple. Indeed, $\phi_n$ are the only non-zero solutions of \eqref{dislo2} vanishing at $x=0$ and $x=1$, so the eigenspaces of $D$ can only be one-dimensional due to the continuity condition at $b$.

\begin{prop} \label{casebhalf}
Let $b=\frac12$. Then, $\operatorname{Spec}(D)=\{\lambda_n\}_{n\in \mathbb{Z}}\subset \mathbb{R}$ where
\begin{equation} \label{ei+half}
\left(\cosh \frac{k}{2}\right)\left(\sin\frac{k}{2}\right)=\left(\cos \frac{k}{2}\right)\left(\sinh\frac{k}{2}\right).
\end{equation} Moreover, for $n\in \mathbb{N}$, $\lambda_{-n}=-\lambda_n$  and
\begin{equation} \label{e+ashalf}
\begin{gathered}
k_n=2\pi\Big(n+\frac14\Big)+\gamma_n \qquad \text{for} \\ \gamma_n=O\left(\exp\left[-2\pi n\right]\right) \qquad \text{as} \quad n\to +\infty.
\end{gathered}
\end{equation}
\end{prop}
\begin{proof} Let us first show that $\lambda=0$ is an eigenvalue. Indeed, when $b=\frac{1}{2}$, the piecewise linear function \eqref{p-wlinear} is such that $\phi\in\operatorname{Dom}(D)$ if and only if $\phi(x)=\alpha \phi_0(x)$ for
\begin{equation} \label{symetricefu}\phi_0(x)=\begin{cases} x & 0\leq x \leq \frac12 \\ 1-x  & \frac12 \leq x \leq 1. \end{cases}\end{equation}
This shows that $0\in \operatorname{Spec}(D)$ is a simple eigenvalue with eigenfunction $\phi_0$.

Now, assume $\lambda\not=0$. The function $\phi(x)$ in \eqref{p-wtrig} is an eigenfunction of $D$ if and only if \eqref{dislo1} holds. Then, $k\in\mathbb{C}\setminus \{0\}$ is such that
\[\begin{bmatrix} \sin \frac{k}{2} & -\sinh\frac{k}{2} \\ k\cos\frac{k}{2} & -k\cosh\frac{k}{2} \end{bmatrix} \begin{bmatrix}\alpha \\ \beta\end{bmatrix} = \begin{bmatrix} 0 \\ 0 \end{bmatrix}
\quad \text{for some} \quad \alpha^2+\beta^2\not=0.\]
This is equivalent to \eqref{ei+half} 
and confirms the first claim.

For the second statement, observe that for $k$ an odd multiples of $\pi/2$, the cosine on the right hand side vanishes, so \eqref{ei+half} is not satisfied. Write the equation in the form
\[
    \tan \frac{k}{2}=\tanh \frac{k}{2}
\]
and assume without loss of generality that $k>0$, as both sides are odd functions. Compare this equation with \eqref{RHS}. Then, note that $\tan\frac{k}{2}=1$ for $k>0$ if and only if $k=2\pi(n+\frac14)$ for some $n\in\mathbb{N}\cup\{0
\}$. Moreover, as $k\to +\infty$,
\[\tanh\frac{k}{2}=1-\frac{2\exp[-\frac{k}{2}]}{\exp[\frac{k}{2}]+\exp[-\frac{k}{2}]}=1+O\left(\exp[-k] \right).\]
Thus,
\[
k_n=2\pi\left(n+\frac14\right) +O\left(\exp\left[ -2\pi n \right]\right) \qquad \qquad n\to+\infty.
 \]
This gives \eqref{e+ashalf}.
\end{proof}

For $c\in(0,1)$, here and elsewhere below we will write
\begin{equation} \label{theepsilon}
     \nu_n\equiv \nu_n^c=\frac{\pi}{c}\Big(n+\frac{\operatorname{sgn}(n)}{4}\Big).
\end{equation}
For our purposes, either $c=b$ or $c=1-b$.
We will suppress the variable $c$ when the context makes it sufficiently clear, mainly in the proofs of the statements.

\begin{prop} \label{prop3} Let $b\in(0,1)$ be such that $b\not=\frac12$. Then, \[\operatorname{Spec}(D)=\{\lambda_n\}_{n\in \mathbb{Z}\setminus\{0\}},\] where
\[\lambda_n=\begin{cases} [k_n^{b}]^2 & n>0 \\ -[k_n^{1-b}]^2 & n<0 \end{cases}\]
and $k_n^{c}$ are the non-zero roots of
\begin{equation} \label{ei+}
\cosh(k(1-c))\sin (kc)=\cos(kc)\sinh(k(1-c)).
\end{equation}
Moreover,
\begin{equation} \label{e+as}\begin{gathered}
k_n^{c}=\nu_n^c+\gamma_n^{c} \qquad \text{where}\\
 \gamma_n^{c}= O\left(\exp\left[-2\pi\frac{1-c}{c} |n|\right]\right) \qquad |n|\to \infty
\end{gathered}\end{equation}
\end{prop}
\begin{proof}
We start with the first claim of this proposition. The piecewise linear function \eqref{dislo2} does not satisfy both conditions in \eqref{dislo1} simultaneously when $b\not=\frac12$. Therefore, for $b\not=\frac12$, $\lambda=0$ is not an eigenvalue of $A$.

The solutions $\phi(x)$ of \eqref{dislo2}  should satisfy \eqref{dislo1}. For $\lambda=k^2>0$, we have
\[\begin{bmatrix} \sin (kb) & -\sinh(k(1-b)) \\ k\cos(kb) & -k\cosh(k(1-b)) \end{bmatrix} \begin{bmatrix}\alpha \\ \beta\end{bmatrix} = \begin{bmatrix} 0 \\ 0 \end{bmatrix}
\quad \text{for} \quad \alpha^2+\beta^2\not=0.\]
For $\lambda=k^2<0$, we have
\[\begin{bmatrix} \sin (k(1-b)) & -\sinh(kb) \\ k\cos(k(1-b)) & -k\cosh(kb) \end{bmatrix} \begin{bmatrix}\alpha \\ \beta\end{bmatrix} = \begin{bmatrix} 0 \\ 0 \end{bmatrix}
\quad \text{for} \quad \alpha^2+\beta^2\not=0.\]
Equating the determinants of these matrices to zero, reduces to \eqref{ei+} with $c=b$ or $c=1-b$, respectively. This confirms the first claim.

For the second claim, assume without loss of generality that $k>0$ and first note that, if $\cos(kc)=0$, the equation is not satisfied. So, \eqref{ei+} can be written as
\[
   \tan(kc)=\tanh(k(1-c)).
\]
Note that $\tan(kc)=1$ for $k>0$ if and only if $k=\frac{\pi}{c}(n+\frac14)$ for some $n\geq 0$. Moreover, as $k\to +\infty$,
\[\tanh(k(1-c))=1-\frac{2\exp[-k(1-c)]}{\exp[k(1-c)]+\exp[-k(1-c)]}=1+O\left(\exp[-2k(1-c)] \right).\]
This gives \eqref{e+as} for $n>0$.
\end{proof}

We will use the following convention for eigenfunctions of $D$. Let \[ B_n^{c}= \frac{\sin(k_n^{c}c)}{\sinh(k_n^{c}(1-c))}.\]
For $n>0$, set
\[ \phi_n(x)=\begin{cases} \sin\big(k_n^{b} x\big) & 0\leq x \leq b \\ B_n^{b} \sinh\big(k_n^{b}(1-x)\big) & b\leq x \leq 1,       \end{cases} \]
where $k_n^{b}$ is the $n$-th positive root of \eqref{ei+} with $c=b$.
For $n<0$, set
\[ \phi_n(x)=\begin{cases} B_n^{1-b} \sinh\big(k_n^{1-b} x\big) & 0\leq x \leq b \\  \sin\big(k_n^{1-b}(1-x)\big) & b\leq x \leq 1,\end{cases}          \]
where $k_n^{1-b}$ is the $n$-th negative root of \eqref{ei+} with $c=1-b$. By construction, $\phi_n$ satisfy \eqref{dislo2}--\eqref{dislo1}, therefore $\phi_n\in\operatorname{Dom}(D)\setminus\{0\}$ and they are the eigenfunctions of $D$. For $b=\frac12$, $\phi_0(x)$ has already been defined in \eqref{symetricefu}.
For notational convenience, in the case $b\not=\frac12$, we define $\lVert\phi_0\rVert^{-1}\phi_0(x)=0$.

Now consider the asymptotic behaviour of
the eigenfunctions of  $D$. 
\begin{prop} \label{prop4}
The eigenfunctions of the operator $D$ can be decomposed as
\[
    \phi_n(x)=\psi_n(x)+\rho_n(x),
\]
where the continuous functions $\psi_n$ and $\rho_n$ are as follows.
For $n>0$,
\begin{equation}
    \label{Positive Eigenfunctions 1}
    \psi_{n} (x) = \begin{cases} \sin\big(\nu_{n}^{b}x\big) \quad  &\qquad x\in[0,b] \\
    \frac{\sqrt{2}}{2}(-1)^{n} \exp\big(-\nu_{n}^{b}(x-b)\big) &\qquad x\in[b,1],
    \end{cases}
\end{equation}
and
\begin{equation} \label{asymprhon}
\|\rho_{n}\|_{\infty}=O( \mathrm{e}^{-\pi\frac{1-b}{b}n}) \qquad n\to\infty.
\end{equation}
For $n<0$,
\begin{equation} \tag{\ref{Positive Eigenfunctions 1}'}
    \label{Negative Eigenfunctions 1}
    \psi_{n} (x) = \begin{cases}
    \frac{\sqrt{2}}{2}(-1)^{n} \exp\big(-\nu^{1-b}_n(b-x)\big) &\qquad x\in[0,b] \\
    \sin\big(\nu_{n}^{1-b}(1-x)\big)  &\qquad x\in[b,1],
    \end{cases}
\end{equation}
and
\begin{equation} \label{asymprhone} \tag{\ref{asymprhon}'}
\|\rho_{n}\|_{\infty}=O( \mathrm{e}^{-\pi\frac{b}{1-b}|n|}) \qquad n\to-\infty.
\end{equation}
Moreover
\begin{equation}
\label{Positive eignef norms}
\|\phi_{n}\|^{-2} = \begin{cases}
\frac{2}{b} + a_{n} & n \to \infty \\
\frac{2}{1-b} + a_{n} &  n\to -\infty,
\end{cases}
\end{equation}
where $a_{n} = O(|n|^{-1})$.
\end{prop}
\begin{proof}
We only include the proof for the case $n>0$, as the other case has a similar proof. 

First note that \[\sin\big(\nu_{n}^{b}b\big)=\sin\big(\nu_{n}^{1-b}(1-b)\big)=\frac{\sqrt{2}}{2}(-1)^{n}\]
hence that $\psi_n\in \Czo$.

Let $x\in[0,b]$.
Then
\begin{align*}
  \phi_n(x)&=\sin(k_n x)=\sin((\nu_n+\gamma_n)x) \\
  &=\sin(\nu_n x) \cos(\gamma_n x)+\sin(\gamma_n x)\cos(\nu_n x) \\
  &= \sin(\nu_n x) +\rho_n (x)
\end{align*}
where
\[
\rho_n(x)=\sin(\nu_n x)(\cos (\gamma_n x)-1)+\sin(\gamma_n x)\cos(\nu_n x).
\]
Recall the asymptotic \eqref{e+as}.
Since, for fixed $x\in [0,b]$, we have
\[
    |\sin(\gamma_n x)|\leq 2|\gamma_n |= O\Big(\exp\Big[-2\pi \frac{1-b}{b} n\Big]\Big),
\]
and the second summing term in the expression for $\rho_n(x)$ is the leading order, we obtain the first claim.

Now, let $x\in [b,1)$. Then,
\begin{align*}
    \phi_n(x)
    &=\sin(k_n b)\frac{e^{k_n(1-x)}-e^{-k_n(1-x)}}{e^{k_n(1-b)}-e^{-k_n(1-b)}} \\
    &= \sin(k_n b)\frac{e^{-k_n (x-b)}-e^{-k_n(2-x-b)}}{1-e^{-k_n(2-2b)}} \\
    &= \left(\sin(\nu_n b)+\bigoh{e^{-2\pi\frac{1-b}{b}n}}\right) \left(e^{-\nu_n (x-b)}+\bigoh{e^{-\pi\frac{1-b}{b}n}}\right) \\ &=\pm \frac{\sqrt{2}}{2} e^{-\nu_n (x-b)}+O(e^{-\frac{1-b}{b}\pi n}).
\end{align*}
This gives the second claim and confirms \eqref{Positive Eigenfunctions 1}. The proof of \eqref{Positive eignef norms} follows from taking the integral square in \eqref{Positive Eigenfunctions 1}. 
\end{proof}

\medskip

We are now ready to complete the proof of Theorem~\ref{Revival in (0,b)}. 

Consider the expansion of the solution of \eqref{disPDE} as a series in the basis of eigenfunctions of $D$, namely
\begin{equation}
    \label{Dislocated solution}
    u(x,t) = \sum_{n\in\mathbb{Z}} \frac{\langle u_{0}, \phi_{n} \rangle}{\|\phi_{n}\|^{2}} e^{-i\la_{n} t} \phi_{n}(x),\qquad t\in \mathbb{R}.
\end{equation}
Recall the convention  $\frac{1}{\|\phi_0\|}\phi_0(x)=0$ for $b\not=\frac12$.

\begin{proof}[Proof of Theorem~\ref{Revival in (0,b)}]
 Let $x\in(0,b)$ and let $t>0$ be fixed. Consider expression \eqref{U0n}, namely
\[
  \widetilde{u_0}(n)=2\int_0^1
   \left[ u_0(by) \sin\Big( \pi\Big(n+\frac14\Big) y\Big)\right.
   +\left. u_0(b^+)\cos\Big( \pi\Big(n+\frac14\Big) y\Big) \right]\mathrm{d}y.
  \]
We split the proof into 4 steps.
In steps 1-3 we confirm the first claim of the theorem, namely we show that the solution can be decomposed as
\[
u(x,t)=U_{\mathcal{R}}(x,t)+U_{\mathcal{C}}(x,t),
\]
where
\[
U_{\mathcal{R}}(x,t)=\sum_{n=1}^\infty\widetilde{u_0}(n)\re^{-i (\nu_n^b)^2t}\sin( \nu_n^b x)
\]
and $U_{\mathcal{C}}(\argdot,t)\in \ContinuousSpace(0,b)$. 
 In  step~4 we then prove the second part of the theorem.

\underline{Step~1}: we isolate the part of the solution which will be continuous. Start with the eigenvalue expansion \eqref{Dislocated solution}. Recall the notation \eqref{theepsilon}. Then, according to \eqref{Positive Eigenfunctions 1} and \eqref{Negative Eigenfunctions 1} in Proposition~\ref{prop4}, we have
    \begin{align*}
        u(x,t)
        &= \sum_{n=-1}^{-\infty}\frac{\langle u_{0}, \phi_{n} \rangle}{\|\phi_{n}\|^{2}} \mathrm{e}^{-i\lambda_{n} t} \left((-1)^n\frac{\sqrt{2}}{2}\mathrm{e}^{-\nu^{1-b}_n(b-x)} +\rho_n(x)\right) \\ & \qquad \qquad + \sum_{n=1}^{\infty} \frac{\langle u_{0}, \phi_{n} \rangle}{\|\phi_{n}\|^{2}} \mathrm{e}^{-i\lambda_{n} t} \left(\sin(\nu^b_n x)+\rho_n(x)\right)+\frac{\langle u_{0}, \phi_{0} \rangle}{\|\phi_{0}\|^{2}}\phi_0(x) \\
        & = \sum_{n=1}^{\infty} \frac{\langle u_{0}, \phi_{n} \rangle}{\|\phi_{n}\|^{2}} \mathrm{e}^{-i\lambda_n t} \sin(\nu^b_n x) +r_1(x,t),
    \end{align*}
where  $r_1$ is given by
\begin{align*}
   r_1(x,t)&=\frac{\sqrt{2}}{2} \sum_{n=-1}^{-\infty}(-1)^n\frac{\langle u_{0}, \phi_{n} \rangle}{\|\phi_{n}\|^{2}} \mathrm{e}^{-i\lambda_{n} t} \mathrm{e}^{-\nu^{1-b}_n(b-x)}\\ & \qquad \qquad +\sum_{n=-\infty}^{\infty}\frac{\langle u_{0}, \phi_{n} \rangle}{\|\phi_{n}\|^{2}} \mathrm{e}^{-i\lambda_{n} t}\rho_n(x)+\frac{\langle u_{0}, \phi_{0} \rangle}{\|\phi_{0}\|^{2}}\phi_0(x).
\end{align*}
Here, the function $\phi_0$ is continuous and $\rho_n\to 0$ as $|n|\to\infty$ uniformly in $[0,b]$ at an exponential rate. Moreover, the functions $\mathrm{e}^{-\nu^{1-b}_n(b-\, \cdot\, )}$ in the first term also converge to zero exponentially fast, uniformly in any compact subset of $[0,b)$. Then, $r_1(\cdot,t)\in\mathrm{C}(0,b)$ for all $t>0$ by Weierstrass' M-test.

Now we need to split the term $\mathrm{e}^{-i\lambda_{n} t}$ in the expression for $u(x,t)$ above, and isolate the leading contribution. Write $\displaystyle{\lambda_n=[\nu_n^b]^2+2\gamma_n^b\nu_n^b+[\gamma_n^b]^2}$ for $n> 0$. Then,
\begin{equation}\label{partproofth3} u(x,t) = \sum_{n=1}^{\infty}\frac{\langle u_{0}, \phi_{n} \rangle}{\|\phi_{n}\|^{2}} e^{-i[\nu_{n}^{b}]^{2} t} \sin(\nu^b_n x)+ r_1(x,t) +r_2(x,t), \end{equation}
where
\begin{align*}
   r_2&(x,t) =
   -i \sum_{n=1}^{\infty} (2\gamma_{n}^b\nu_n^b+[\gamma_n^b]^2) \frac{\langle u_{0}, \phi_{n} \rangle}{\|\phi_{n}\|^{2}} e^{-i[\nu_{n}^{b}]^{2} t} \left[ \int_{0}^{t}e^{-i(2\gamma_{n}^b\nu_n^b+[\gamma_n^b]^2) s} \mathrm{d}s\right] \sin(\nu^b_n x).
\end{align*}
By virtue of \eqref{e+as}, \eqref{Positive eignef norms}
and Weierstrass' test once again, $r_2$ is also a continuous function of $x\in(0,b)$.

Finally, consider the series on the right hand side of \eqref{partproofth3}. The generalised Fourier coefficients are given by
\begin{align*}
    \frac{\langle u_{0}, \phi_{n} \rangle}{\|\phi_{n}\|^{2}}&=(\alpha_n+\beta_n+\langle u_0,\rho_n\rangle)\left(\frac{2}{b}+a_n\right)
    =\frac{2}{b}(\alpha_n+\beta_n)+\delta_n
\end{align*}
where $a_n=\bigoh{n^{-1}}$ from \eqref{Positive eignef norms},
\begin{gather}
\label{alphan}
\begin{aligned}\frac{2}{b} \alpha_n&=\frac 2 b\int_0^b u_0(y)\sin(\nu_n^b y)\,\mathrm{d}y \\ &=2\int_0^1u_0(by)\sin\Big( \pi\Big(n+\frac14\Big) y\Big)\mathrm{d}y, \end{aligned}  \\ \nonumber
\frac 2 b\beta_n=\frac{\sqrt{2}}{b}(-1)^n\int_b^1 u_0(y)\mathrm{e}^{-\nu_n^b(y-b)}\,\mathrm{d}y \qquad \text{and} \\ \nonumber \delta_n=\frac{2}{b}\langle u_0,\rho_n\rangle+(\alpha_n+\beta_n+\langle u_0,\rho_n\rangle)a_n.
\end{gather}
We now write
\begin{align*}
     \frac{2}{b} \sum_{n=1}^{\infty}  \beta_n &e^{-i[\nu_{n}^{b}]^{2} t} \sin(\nu^b_n x)=
    r_3(x,t)+ 
    \\& \frac{2u_0(b^+)}{b}\sum_{n=1}^\infty \left[\int_0^b \cos\left(\frac{\pi}{b} \Big(n+\frac14\Big) y\right)\mathrm{d}y\right]
     \mathrm{e}^{-i\left[\frac{\pi(n+\frac14)}{b}\right]^2 t} \sin\left(\frac{\pi}{b} \Big(n+\frac14\Big) x\right)
\end{align*}
and
\[
    \sum_{n=1}^{\infty}  \delta_n \mathrm{e}^{-i[\nu_{n}^{b}]^{2} t} \sin(\nu^b_n x)=r_4(x,t).
\]

\smallskip

\underline{Step~2}: we prove that $r_3$ and $r_4$ are continuous functions of $x\in(0,b)$. We first show the following asymptotic formulas, valid as $n\to \infty$:
\begin{gather} \label{claimasymp} \alpha_n=\bigoh{n^{-1}}, \quad
    \beta_n=\frac{(-1)^n b}{\sqrt{2}\pi} u_0(b^+)\Big(n+\frac14\Big)^{-1}+\bigoh{n^{-2}}  \\ \nonumber \text{and}  \qquad \delta_n=\bigoh{n^{-2}}.
\end{gather}
    Note that the last formula is a consequence of the first and second one, together with the  facts that $a_n=\bigoh{n^{-1}}$ and that, by Proposition~\ref{prop4}, $\rho_n\to 0$ as $|n|\to\infty$ uniformly at an exponential rate. 
    
In order to show the two identities in \eqref{claimasymp}, we first decompose $u_0$ into an absolutely continuous function $u^{\mathrm{ac}}_0$ and a piecewise constant function $u_0^j$ with jumps at finitely many point $\{x_1,...,x_{k}\}\subset (0,1)$, see Remark~\ref{crucial}.  Then we have
\begin{align*}
   \alpha_n&=-\frac{1}{\nu_n^b}\int_0^b u_0(y)\,\,\mathrm{d}\!\!\left[\cos(\nu_n^b y)\right] \\ &=-\frac{1}{\nu_n^b}\left[(-1)^n u_0(b^-)\frac{\sqrt{2}}{2}-u_0(0^+)-\int_0^b\cos(\nu_n^b y)\,\mathrm{d}u_0(y)\right],
\end{align*}
where $\mathrm{d}u_0(y)$ is the Lebesgue-Stieltjes measures associated to $u_0$ in the sub-interval $[0,b)$, so that
$$
\int_0^b\cos(\nu_n^b y)\,\mathrm{d}u_0(y)=\int_0^b\cos(\nu_n^b y)\,\mathrm{d}u^{\mathrm{ac}}_0(y)+\sum_{j=1}^{m-1}h_j\cos(\nu_n^b x_j)
$$
where $\{x_1,\ldots,x_m-1\}\subset (0,b)$ are the points of discontinuity in this sub-interval, with jumps $h_j$.
It follows from the boundedness of $\mathrm{d}u^{\mathrm{ac}}_0$ that
\[
\left|\int_0^b\cos(\nu_n^b y)\mathrm{d}u_0(y)\right|
<\infty.
\]
Then, since $\nu_n^b=O(n)$, indeed we have $\alpha_n=O(n^{-1})$ as $n\to\infty$, proving the first asymptotic statement  in \eqref{claimasymp}.

Let us now consider the asymptotic behaviour of $\beta_n$.
    Since \[\psi_n(y)\,\mathrm{d}y=-\frac{1}{\nu_n^b}\mathrm{d}\psi_n(y)\] for all $y\in[b,1]$, we have
    \begin{align*}
        \beta_n&=\frac{\sqrt{2}}{2} (-1)^n \frac{-1}{\nu_n^b}\int_b^1 u_0(y)\,\mathrm{d}\psi_n(y) \\ &=
        \frac{\sqrt{2}}{2} (-1)^n\frac{1}{\nu_n^b}\left[u_0(b^+)\mathrm{e}^{-\nu_n^b(1-b)}+\int_b^1\mathrm{e}^{-\nu_n^b(y-b)} \mathrm{d}u_0(y) \right]
    \end{align*}
    where now 
$\mathrm{d}u_0(y)$ is the Lebesgue-Stieltjes measure associated to $u_0$ in the sub-interval $(b,1]$, so that we have (see Remark \ref{crucial}) 
    \[
    \int_b^1\mathrm{e}^{-\nu_n^b(y-b)} \mathrm{d}u_0(y)=\int_b^1\mathrm{e}^{-\nu_n^b(y-b)} \mathrm{d}u^{\mathrm{ac}}_0(y) +\sum_{j=m}^k h_j\mathrm{e}^{-\nu_n^b(x_j-b)}
    \]
    where $\{x_m,\ldots,x_k\}\subset (b,1)$ are the points of discontinuity with jumps $h_j$. The finite summation decays exponentially and
    \begin{align*}
       \left| \int_b^1\mathrm{e}^{-\nu_n^b(y-b)} \mathrm{d}u^{\mathrm{ac}}_0(y)\right|& \leq \int_{b}^1\mathrm{e}^{-\nu_{n}^b(y-b)} |(u^{\mathrm{ac}}_0)'(y)| \,\mathrm{d}y \\
       &\leq \| (u^{\mathrm{ac}}_0)'\|_{\infty}
        \int_{b}^1\mathrm{e}^{-\nu_{n}^b(y-b)} \,\mathrm{d}y = \bigoh{n^{-1}} 
    \end{align*}
    as $n\to \infty$. Hence, indeed $\beta_n$ has the asymptotic behaviour claimed in \eqref{claimasymp} and this completes the proof of the latter. 
    
    Thus, since the reminders of $\beta_n$ and $\delta_n$ are $\bigoh{n^{-2}}$, then indeed $r_3$ and $r_4$ are  continuous functions.

\smallskip

\underline{Step~3}: we now complete the proof of the first part of the theorem. The leading order of $\frac 2 b \beta_n$  is given by\be\label{betan}
\frac{2}{\pi}u_0(b^+)\frac{(-1)^n}{\left(n+\frac14\right)} \frac{\sqrt{2}}{2}
=2 u_0(b^+)\int_0^1 \cos \left[\left(n+\frac14\right)\pi y\right]\mathrm{d}y.
\ee
Hence we have
\begin{align*}
\frac{2}{b}\sum_{n=1}^\infty &\left[\frac{(-1)^n b\sqrt{2}}{2\pi\left(n+\frac14\right)}u_0(b^+)\right] \mathrm{e}^{-i(\nu_n^b)^2t}\sin(\nu_n^b x)
 \\ &=
\frac{2}{\pi}u_0(b^+)\sum_{n=1}^\infty (-1)^n \frac{\sqrt{2}}{2}\mathrm{e}^{-i(\nu_n^b)^2t}\frac{\sin(\nu_n^b x)}{n+\frac14}
\\
&= 2 u_0(b^+) \sum_{n=1}^\infty\left(\int_0^1 \cos \left[\left(n+\frac14\right)\pi y\right]\mathrm{d}y\right)
\mathrm{e}^{-i(\nu_n^b)^2t}\sin(\nu_n^b x).
\end{align*}
Considering the expression for $\alpha_n$ given in \eqref{alphan} and the contribution from $\beta_n$ in \eqref{betan}, the above implies the expression \eqref{U0n} for $\widetilde{u_0}(n)$. Therefore, the first claim of Theorem~\ref{Revival in (0,b)} is valid.

\medskip

\underline{Step~4}: we prove the second part of the theorem.

We aim to expressf $\widetilde{u_0}(n)$ in terms of the Fourier coefficient of a single function, defined on the double interval $(0,2)$. Expanding sine and cosine terms, and changing variables to $z=2-y$, we obtain
\begin{align*}
 \widetilde{u_0}(n) &= \frac 1 2\int_0^1 [u_0(b^+)+iu_0(by)] \mathrm{e}^{-i\pi (n+\frac14) y}\mathrm{d}y+
      \frac 1 2 \int_0^1 [u_0(b^+)-iu_0(by)]\mathrm{e}^{i\pi (n+\frac14) y}\mathrm{d}y \\
       &=\frac 1 2\int_0^1 [u_0(b^+)+iu_0(by)] \mathrm{e}^{-i\frac{\pi}{4} y}\mathrm{e}^{-i\pi n y}\mathrm{d}y+ \\
       & \qquad \qquad \qquad \frac 1 2\int_1^2 [u_0(b^+)-iu_0(2b-bz)]\mathrm{e}^{i\pi (n+\frac14) (2-z)}\mathrm{d}z \\
       &=\frac 1 2\int_0^1 [u_0(b^+)+iu_0(by)] \mathrm{e}^{-i\frac{\pi}{4} y}\mathrm{e}^{-i\pi n y}\mathrm{d}y+ \\
       & \qquad \qquad \qquad \frac 1 2\int_1^2 [u_0(b^+)-iu_0(2b-bz)] \mathrm{e}^{i\frac{\pi}{2}}\mathrm{e}^{-i\frac{\pi}{4} z}\mathrm{e}^{-i\pi n z} \mathrm{d}z \\
       &=\frac 1 2\int_0^2 G_{u_0}^b(y)e^{-i\pi n y}\mathrm{d}y=
       \widehat{G_{u_0}^b}(n)
\end{align*}
where $\widehat \argdot$ denotes the Fourier coefficient on the interval $(0,2)$  and
	\begin{equation}
		\tag{\ref{v,g0s}}
		\begin{aligned}
			&G_{u_0}^b(x) = e^{- i \frac{\pi x}{4}}
		 \begin{cases}
			u_0(b^+)+iu_0(bx), &\quad 0\leq x \leq 1 \\
			iu_0(b^+)+u_0(b(2-x)), &\quad 1\leq x \leq 2.
		\end{cases}
		\end{aligned}
	\end{equation}
	
We now evaluate this expression at rational times. 	For $t=t^{(b)}=\frac{(2b)^2p}{2\pi q}$, we have
\begin{align*}
   U_{\mathcal{R}} (x,t^{(b)})&=\sum_{n=1}^\infty\widetilde{u_0}(n) \mathrm{e}^{-i(\nu_n^b)^2 t^{(b)}}\sin(\nu_n^b x)\\ &=\sum_{n=1}^\infty\widetilde{u_0}(n) \mathrm{e}^{-i(n+\frac14)^2\frac{2\pi p}{q}}\sin\Big(\frac{\pi}{b}\Big(n+\frac14\Big)x\Big).
   \end{align*}
 Using the modularity identity
 \[
    \mathrm{e}^{-i(n+\frac{1}{4})^{2}\frac{2\pi p}{q}} = e^{-i(m+\frac{1}{4})^{2}\frac{2\pi p}{q}}\quad \text{for}\quad n\equiv_{2q} m
    \]
    and the identity \eqref{unitroots} for the sum of all $2q$ roots of unity, after algebraic manipulation analogous to the ones in the proof of lemma~\ref{lem:Airy.MainRevivalLemma}, we find
   \begin{align*}
   U_{\mathcal{R}} (x,t^{(b)}) 
   &=\sum_{m=0}^{2q-1}e^{ -i(m+\frac{1}{4})^{2}\frac{2\pi p}{q} } \frac{1}{2q}\sum_{k=0}^{2q-1}e^{i\pi \frac{mk}{q}}
  \sum_{n=1}^\infty\widetilde{u_0}(n) \mathrm{e}^{-i\frac{\pi k n}{q}}\sin\Big(\frac{\pi}{b}\Big(n+\frac14\Big)x\Big)
  \\
     &=\sum_{k=0}^{2q-1}d_k^{p,q}\sum_{n=1}^\infty\widetilde{u_0}(n) \mathrm{e}^{-i\frac{\pi k n}{q}}\sin\Big(\frac{\pi}{b}\Big(n+\frac14\Big)x\Big),
\end{align*}
where
  \begin{equation}\label{Ckmpq}
   d_{k}^{(p,q)} = \frac{1}{2q}\sum_{m=0}^{2q-1}e^{ -i(m+\frac{1}{4})^{2} \frac{2\pi p}{q}} e^{i\pi \frac{mk}{q}}. \end{equation}
Hence
\begin{equation} \label{finHildis}
\begin{aligned}
\sum_{n=1}^\infty&\widetilde{u_0}(n) \mathrm{e}^{-i\frac{\pi kn}q}\sin\Big(\frac{\pi}{b}\Big(n+\frac14\Big)x\Big) \\ 
& =\sum_{n=1}^\infty \frac{ \widetilde{u_0}(n)}{2i} \left[\mathrm{e}^{i\frac{\pi }{4b}x}\mathrm{e}^{i\pi(\frac{x}{b}-\frac{k}{q})n}- \mathrm{e}^{-i\frac{\pi }{4b}x}\mathrm{e}^{-i\pi(\frac{x}{b}+\frac{k}{q})n}\right]  \\
&=i\mathrm{e}^{-i\frac{\pi }{4b}x}\sum_{n=1}^\infty \frac{ \widetilde{u_0}(n)}{2}\mathrm{e}^{i\pi(-\frac{x}{b}-\frac{k}{q})n}-
i\mathrm{e}^{i\frac{\pi }{4b}x}\sum_{n=1}^\infty \frac{ \widetilde{u_0}(n)}{2}\mathrm{e}^{i\pi(\frac{x}{b}-\frac{k}{q})n}
\end{aligned}
\end{equation}

Since $ \widetilde{u_0}(n)=\widehat{G_{u_0}^b}(n)$, using the canonical identity \eqref{HT identities} for the periodic Hilbert transform on $[0,2]$ gives
\[
\sum_{n=1}^\infty \frac{ \widetilde{u_0}(n)}{2}\mathrm{e}^{i\pi(-\frac{x}{b}-\frac{k}{q})n}=\frac{G_{u_0}^b(-\frac{x}{b}-\frac{k}{q})-\langle G_{u_0}^b\rangle +i[\mathcal{H}G_{u_0}^b](-\frac{x}{b}-\frac{k}{q})}{2}\]
and
\[
\sum_{n=1}^\infty \frac{ \widetilde{u_0}(n)}{2}\mathrm{e}^{i\pi(\frac{x}{b}-\frac{k}{q})n}=\frac{G_{u_0}^b(\frac{x}{b}-\frac{k}{q})-\langle G_{u_0}^b\rangle+i[\mathcal{H}G_{u_0}^b](\frac{x}{b}-\frac{k}{q})}{2}.
\]
Then, \eqref{finHildis} can be rewritten as
\begin{align*}
\sum_{n=1}^\infty & \widetilde{u_0}(n) \mathrm{e}^{-i\frac{\pi kn}q}\sin\Big(\frac{\pi}{b}\Big(n+\frac14\Big)x\Big) =
 i\mathrm{e}^{-i\frac{\pi }{4b}x}\left[\frac{G_{u_0}^b(-\frac{x}{b}-\frac{k}{q})+i[\mathcal{H}G_{u_0}^b](-\frac{x}{b}-\frac{k}{q})}{2}\right]
\\ &\qquad
-
i\mathrm{e}^{i\frac{\pi }{4b}x}\left[\frac{G_{u_0}^b(\frac{x}{b}-\frac{k}{q})+i[\mathcal{H}G_{u_0}^b](\frac{x}{b}-\frac{k}{q})}{2}\right]
 -2\langle G_{u_0}^b\rangle \sin\Big(\frac{\pi x}{4b}\Big).
\end{align*}
Therefore,
$$
   U_{\mathcal{R}} (x,t^{(b)})  =L_1(x)+L_2(x)+L_3(x),
   $$
   where
      \begin{align*}
  L_1(x)   &=\frac 1 2\sum_{k=0}^{2q-1}d_k^{p,q}\left[ i\mathrm{e}^{-i\frac{\pi }{4b}x}G_{u_0}^b\Big(-\frac{x}{b}- \frac{k}{q}\Big)-i\mathrm{e}^{i\frac{\pi }{4b}x}G_{u_0}^b\Big(\frac{x}{b}-\frac{k}{q}\Big)\right],
  \\
 L_2(x)   &=\frac 1 2\sum_{k=0}^{2q-1}d_k^{p,q} \left[ \mathrm{e}^{i\frac{\pi }{4b}x} [\mathcal{H}G_{u_0}^b]\Big(\frac xb -\frac kq\Big) - \mathrm{e}^{-i\frac{\pi }{4b}x} [\mathcal{H}G_{u_0}^b]\Big(-\frac xb -\frac kq\Big)  \right]
 \\
 L_3(x)&=-2\left(\sum_{k=0}^{2q-1}d_k^{p,q}\right) \langle G_{u_0}^b\rangle \sin\Big(\frac{\pi x}{4b}\Big).
  \end{align*}
The simplification
\[
    \sum_{k=0}^{2q-1}d_k^{p,q} = \re^{-i\frac{\pi p}{2q}}
\]  in the expression for $L_3$ completes the proof.
\end{proof}

\begin{remark} \label{reversing}
The reflection operator $Ru_0(x)=u_0(1-x)$ is an isometry of $\Ltwozo$, $R^2=I$ and $R^*=R$. Moreover, in the obvious notation,
\[
    R\phi_{n,b}(x)=\phi_{-n,1-b}(x).
\]
Therefore $R$ is a diffeomorphism between the domains of $D_b$ and $D_{1-b}$.   Since
\begin{align*}
RD_{1-b}Rf&=R\sum_{n\in\mathbb{Z}} \lambda_{n,1-b}\frac{\langle Rf,\phi_{n,1-b}\rangle}{\|\phi_{n,1-b}\|^2} \phi_{n,1-b}\\ &=\sum_{n\in\mathbb{Z}} \lambda_{n,1-b}\frac{\langle f,\phi_{-n,b}\rangle}{\|R\phi_{-n,b}\|^2} \phi_{-n,b} \\
&=-\sum_{n\in\mathbb{Z}} \lambda_{n,b}\frac{\langle f,\phi_{n,b}\rangle}{\|\phi_{n,b}\|^2} \phi_{n,b} \\&=-D_bf
\end{align*}
for all $f\in \mathrm{Dom}(D_b)$, then the operators $D_{1-b}$ and $-D_b$ are similar. Indeed, their spectra map accordingly.
\end{remark}

By virtue of this remark, the solution to \eqref{disPDE} for $x\in(0,b)$ is
\begin{align*}
u(x,t)&=\sum_{n\in\mathbb{Z}}\frac{\langle u_0,\phi_{n,b}\rangle}{\|\phi_{n,b}\|^2} \mathrm{e}^{-i\lambda_{n,b}t}
\phi_{n,b}(x) \\
&=\sum_{n\in\mathbb{Z}}\frac{\langle u_0,R\phi_{-n,1-b}\rangle}{\|R\phi_{-n,1-b}\|^2} \mathrm{e}^{-i\lambda_{n,b}t}
R\phi_{-n,1-b}(x) \\
&=\sum_{n\in \mathbb{Z}} \frac{\langle Ru_0,\phi_{n,1-b}\rangle}{\|\phi_{n,1-b}\|^2} \mathrm{e}^{i\lambda_{n,1-b}t} \phi_{n,1-b}(1-x)\\&=v(1-x,-t)
\end{align*}
where $v$ is the solution to \eqref{disPDE}, with initial condition $Ru_0$, for $x\in(b,1)$. This observation justifies the comment made in Section~\ref{section2} about the formula for the revivals of the solution in the complementary interval $(b,1)$.

\section{Towards a general definition of revivals} \label{rev}
In this final section, we propose a rigorous framework fitting the new revival phenomenon formulated in theorems~\ref{lkdvthm2} and \ref{Revival in (0,b)} within the broader context of revivals in dispersive boundary value problems.

Let $l>0$. Consider a linear dispersive equation of the form
\be\label{gendisp}
u_t+iTu=0,\qquad u(\argdot,t)\in \Lebesgue^2(0,l), \;\; t\in\mathbb{R},
\ee
where $T$ is a linear operator with domain a dense subspace of $\Lebesgue^2(0,l)$, defined via its dispersion relation.  Revivals occur when the time variable is what we normally call a \emph{rational time}. These rational times have the form $t=\frac p q \alpha$, for $p,q\in\mathbb N$ coprime numbers, and $\alpha\in\mathbb R$ a fixed parameter that depends on the structure of the operator $T$ and the length $l$.

The first notion of revivals was established in the 1990s and it describes what has been called {\em Talbot effect} for linear dispersive systems, \cite{berry1996integer}. It is present, for example, in the solutions of the periodic linear Schr\"odinger equation \cite{berry1996integer} and the periodic Airy equation \cite{olver2010dispersive}. This notion prompted most of the recent research on the subject. We propose the following definition to describe this type of revival.

\begin{definition}[The periodic revival property]\label{prp}
A boundary value problem for equation \eqref{gendisp} is said to admit the \emph{periodic revival property}, if
its solution evaluated at any rational time is equal to a finite linear combination of translated and reflected copies of the product of the initial datum with a continuous function.
\end{definition}
In most known cases, the continuous function is a constant.

A different manifestation of revivals,  described in the context of the periodic linear Benjamin-Ono equation in \cite{boulton2020new}, is the following. In this manifestation, initial jump discontinuities give rise to logarithmic cusps. This phenomenon matches the contributions of the term $U_{\mathcal{R}}$ in theorems~\ref{lkdvthm2} and \ref{Revival in (0,b)} at rational times, and it prompts the next definition. 

\begin{definition}[The cusp revival property]\label{prpc}
A boundary value problem for equation \eqref{gendisp} is said to admit the \emph{cusp revival property}, if
its solution evaluated at any rational time is the linear combination two functions. One of these functions is a finite linear combination of translated and reflected copies of the product of a continuous function with the initial datum. The other function is a finite linear combination of translated and reflected copies of the Hilbert transform of the product of a continuous function with the initial datum.
\end{definition}

Our third definition describes a revival property modulo a continuous contribution to the solution at rational times.

\begin{definition}[The weak revival property]\label{wprp}
A boundary value problem for equation \eqref{gendisp} is said to admit the \emph{weak revival property}, if its solution is the sum of a function possessing a revival property of the form described in the definitions~\ref{prp} or \ref{prpc}, plus a continuous function.
\end{definition}

This weaker form of revival fits the phenomena described above for the boundary value problems \eqref{pseD1} and \eqref{disPDE}. A manifestation of it was first reported in \cite{rodnianski1999continued} and then examined in the  case of the linear Schr\"odinger equation with specific Robin-type boundary conditions in \cite{BFP}.
Moreover, this effect has been found in several other linear and nonlinear equations that are natural generalisations of the periodic problems originally studied by Berry in \cite{berry1996quantum,berry1996integer} and by Olver in \cite{olver2010dispersive}.

\medskip

Below we list currently known boundary value problems that exhibit the revival effects described by the above definitions. Some of these involve results formulated only for specific discontinuous initial profiles. 
\begin{itemize}
\item \emph{Periodic revival.} Linear constant coefficient dispersive equations with periodic boundary conditions. Linear Schr\"odinger equation with pseudo-periodic, Dirichlet or Neumann boundary conditions. 
\item \emph{Weak periodic revival.}
The linear Schr\"odinger equation with certain Robin boundary conditions \cite{BFP}, with additional numerical evidence for the case of general Robin boundary conditions \cite{olver2018revivals}. 
The linear Schr\"{o}dinger equation with a complex potential, on the 1D-torus \cite{rodnianski1999continued, cho2021talbot} and with Dirichlet boundary conditions \cite{boulton2023phenomenon}.
The periodic problem for bidirectional dispersive hyperbolic equations such as the linear beam equation \cite{farmakis2023new}.
The periodic problem for the nonlinear Schr\"odinger  and the Korteweg-deVries   equations \cite{erdogan2013talbotpaper,erdougan2013global}. Numerical evidence for strongly nonlinear, non-integrable generalisations of the Korteweg-deVries equation \cite{chen2013dispersion}.
\item \emph{Cusp revival.} Only the linear Benjamin-Ono equation \cite{boulton2020new}.
\item \emph{Weak cusp revival.} 
The periodic linearised Intermediate Long Wave and Smith equations \cite{boulton2020new}. Numerical evidence for finite time for the periodic nonlinear Benjamin-Ono equation \cite{chen2013dispersion}. 
\end{itemize}

\appendix

   
\section{Notation} \label{apa}

In this appendix we give details about the notation that we used in this paper. 

\subsubsection*{Piecewise Lipschitz functions}
In this work a function $u_0$ is called piecewise Lipschitz, if there exist \[\{x_0=0<x_1<\cdots<x_k<x_{k+1}=1\}\subset [0,1]\] such that  $u_0$ is Lipschitz in the open interval $(x_{j},x_{j+1})$ for all $j=0,\ldots,k$. This ensures in particular that $u_0\in \operatorname{BV}([0,1])$ and $u_0'\in \operatorname{L}^{\infty} (x_j,x_{j+1})$.

\begin{remark} \label{crucial}
If $u_0$ is piecewise Lipschitz, then $u_0=u^{\mathrm{ac}}_0+u^{\mathrm{j}}_0$ 
where $u^{\mathrm{ac}}_0\in \mathrm{AC}([0,1])$, $(u^{\mathrm{ac}}_0)'\in \mathrm{L}^{\infty}(0,1)$ and 
 $u^{\mathrm{j}}_0$ is the finite linear combination of the characteristic functions of intervals. 
 In particular, regarded as a function of bounded variation in $[0,1]$, $u_0$ has only finitely many jump discontinuities and the singular part in its Lebesgue-Stieltjes decomposition is zero.
\end{remark}

\subsubsection*{Fourier series}
For $l>0$, let $u\in\Lebesgue^2(0,l)$. The Fourier coefficients of $u$ are denoted by
\[
\widehat{u}(n) = \frac{1}{l} \int_{0}^{l} u(y) e^{-i\frac{2\pi}l ny} \mathrm{d}y, \qquad \qquad n\in\mathbb{Z}.
\]
Writing the mean of the function $u$ as $\langle u\rangle=\widehat{u}(0)$, we have
\[u(x) = \langle u\rangle + \sum_{n=1}^{\infty} \Big[\widehat{u}(-n) e^{-i\frac{2\pi}{l}nx} + \widehat{u}(n) e^{i\frac{2\pi}{l}nx}\Big].\]

\subsubsection*{Periodic Hilbert transform}
The periodic Hilbert transform on $(0,l)$ denoted by
$\mathcal{H}:\Lebesgue^2(0,l)\longrightarrow \Lebesgue^2(0,l)$ is given by
        \begin{equation}
            \label{General HT}
            \mathcal{H}u(x) = i\sum_{n=1}^{\infty}\Big[ \widehat{u}(-n) e^{-i\frac{2\pi}{l}nx} - \widehat{u}(n) e^{i\frac{2\pi}{l}nx}\Big].
        \end{equation}
       See \cite{king} and references therein.
       
The canonical identity,
\begin{equation} \label{HT identities}
\sum_{n=1}^{\infty} \widehat{u}(n) \mathrm{e}^{i\frac{2\pi}{l}nx} = \frac{u(x)-\langle u\rangle +i\mathcal{H}u(x)}{2},
\end{equation}
crucial in the proofs of theorems~\ref{Revival in (0,b)}   and \ref{lkdvthm2}, is a routine consequence of this definition.

The Hilbert transform of a function of bounded variation with only a finite number of jump discontinuities and no singular part, displays a logarithmic cusp singularity at any point where the function has a jump discontinuity, see the illustration in Figure~\ref{fig0}. 
This follows from the Lebesgue Decomposition Theorem  \cite[Section 33.3]{kolmogorov1975introductory}, the linearity of $\mathcal{H}$ and the fact that, for $0<a<b<l$,
\begin{equation} \label{logsing}
    \mathcal{H} \mathbbm{1}_{[a,b]}(x)=\frac{l}{2}\log\left|\frac{\sin\Big(\pi\frac{x-a}{l}\Big)}{\sin\Big(\pi\frac{x-b}{l}\Big)}\right|.
\end{equation}
The latter is shown by a direct calculation, using the principal value representation of $\mathcal{H}$ as a periodic Calderon-Zygmund operator,
\be\label{HilbOp}
 \mathcal{H}u(x) =\frac 1 {l}
 \dashint_{0}^l u(y)\cot\left(\frac{\pi(x-y)}{l}\right)\, \mathrm{d}y.
\ee


\section{The spatial operators} \label{apb}
The existence and uniqueness of solutions for the boundary value problems considered in this paper is a direct consequence of the self-ad\-joint\-ness and compactness of the resolvents of the spatial operators. In this appendix we give details about the proof of these two properties. 

\subsubsection*{The spatial operator in \eqref{pseD1}}
Let $A:\mathrm{Dom}(A)\longrightarrow \mathrm{L}^2(0,1)$ be given by the differential expression
$$
A\phi(x)=i\phi_{xxx}(x)
$$
on the domain defined by \eqref{airydom1}, as described in Section~\ref{sec4}. Then, $A$ is self-adjoint and hence $iA$ is the generator of a $C_0$ one-parameter semigroup. Moreover, its resolvent is compact. Therefore, the solutions of \eqref{pseD1} are given as a series in the orthonormal basis of eigenfunctions of $A$ for any $u_0\in  \mathrm{L}^2(0,1)$.

We prove that $A$ is self-adjoint as follows.  Firstly note that the minimal operator associated to the differential expression  $i\partial_{x}^3$ is symmetric and has equal deficiency indices with value 3. The boundary form associated to this minimal operator is
\[
    \{\phi,\psi\}=[\phi''(x)\overline{\psi(x)}-\phi'(x)\overline{\psi'(x)}+\phi(x)\overline{\psi''(x)}]\Big|_{x=0}^{x=1}.
\]
Since the boundary conditions in \eqref{airydom1} are linearly independent and $\{\phi,\psi\}=0$ for any $\phi,\psi\in \operatorname{Dom}(A)$, they constitute a symmetric set of boundary conditions for the minimal operator. Hence, since the operator $A$ is the restriction of the maximal operator (the adjoint of the minimal operator) to the subspace generated by these boundary conditions, it indeed follows \cite[Theorem~XII.4.30]{DunfSchv2} that $A$ is a self-adjoint operator. Thus, according to \cite[Theorem~XIII.4.1]{DunfSchv2}, $A$ has a compact resolvent. 

\medskip
\subsubsection*{The spatial operator in \eqref{disPDE}}
Let 
\[
     \mathrm{Dom}(D)=\Big\{u\in \operatorname{AC}([0,1]): \operatorname{sgn}(\cdot-b)u'\in \operatorname{AC}([0,1]),\,
u(0)=u(1)=0\Big\}.
\]
We claim that the dislocated Dirichlet Laplacian $D:\mathrm{Dom}(D)\longrightarrow \mathrm{L}^2(0,1)$ given by the quasi-differential expression \eqref{expressiond}, is also self-adjoint. The proof of this is less standard.
First, note that $\mathrm{Dom}(D)$ is a dense subspace of $\Ltwozo$, since the restriction of functions in $\operatorname{Dom}(D)$ to the sub-intervals $(0,b)$ and $(b,1)$, are dense in $\Lebesgue^2(0,b)$ and $\Lebesgue^2(b,1)$ respectively. The operator $D$ fits into the classical framework of quasi-differential operators. Since $\operatorname{sgn}(x-b)$ only vanishes at $x=b$, the expression in \eqref{expressiond} is regular in $[0,1]$, see \cite[\S III.10.2]{EdmundsEvans}. Moreover,  the associated minimal operator is symmetric and therefore it forms a compatible pair with itself. As the operator $D$ is regularly solvable with respect to this pair, then indeed $D=D^*$ by virtue of  \cite[Theorem~III.10.6]{EdmundsEvans}.

We conclude this appendix by giving a full proof that $D$ has a compact resolvent.

\begin{prop}
    The resolvent of $D$ is compact.
\end{prop}
\begin{proof}
Since $D$ is self-adjoint, $i$ is in its resolvent. We aim to compute the integral kernel associated to the inverse of $D-i$, via variation of parameters, then show that this operator is compact. Let $\omega=\mathrm{e}^{i\frac{\pi}{4}}$, so that $\omega^2=i$. Let
\[
\Phi_1(x)=\begin{cases} \mathrm{e}^{i\frac{3\pi}{8}}\cos(\omega (x-b)) & 0\leq x <b \\
 \mathrm{e}^{i\frac{3\pi}{8}}\cosh(\omega (x-b)) & b< x \leq 1 \end{cases}
\]
and
\[
\Phi_2(x)=\begin{cases} \mathrm{e}^{i\frac{3\pi}{8}}\sin(\omega (x-b)) & 0\leq x <b \\
 -\mathrm{e}^{i\frac{3\pi}{8}}\sinh(\omega (x-b)) & b< x \leq 1. \end{cases}
\]
Then, $(\text{sgn}(x-b) \Phi_j'(x))'=i\Phi_j(x)$ and the Wronskian
\[
  \begin{aligned}   W_0(\Phi_1,\Phi_2)&=\Phi_1(0)\sgn(0-b)\phi_2'(0)-\sgn(0-b)\phi_1'(0)\phi_2(0) \\ &=\frac{-1}{\omega} (-1)[\omega \cos^2(\omega b)+\omega \sin^2(\omega b))]=1. \end{aligned}
\]
Note that $\Phi_j$ are continuous at $b$ and $\Phi_j'(b^-)=-\Phi_j'(b^+)$.

The variation of parameters formula, valid for regular quasi-differential operators \cite[Lemma~III.10.9]{EdmundsEvans}, gives that $(\text{sgn}(x-b) u'(x))'-i u(x)=f(x)$ for $f\in \mathrm{L}^1_{\mathrm{loc}}(0,1)$, if and only if
\begin{equation} \label{varpa}
\begin{aligned}
    u(x)&=c_1(f)\Phi_1(x)+c_2(f)\Phi_2(x)\\ &\qquad +\int_0^x [\Phi_1(x)\Phi_2(y)-\Phi_2(x)\Phi_1(y)]f(y)\,\mathrm{d}y
\end{aligned}\end{equation}
where $c_j(f)\in \mathbb{C}$ are suitable constants, depending on $f$. The expression ensures that $u(\argdot)$ and
$\sgn(\argdot-b)u'(\argdot)$ are absolutely continuous in $[0,1]$, where  $u(b^-)=u(b^+)$ and $u'(b^-)=-u'(b^+)$. We now give the unique $c_j(f)$, so that $u(0)=u(1)=0$ for $f\in \mathrm{L}^2(0,1)$. This will complete the proof.

Let
\[
    M=\begin{bmatrix} \Phi_1(0) & \Phi_2(0) \\ \Phi_1(1) & \Phi_2(1) \end{bmatrix}=
    \mathrm{e}^{i\frac{3\pi}{8}} \begin{bmatrix} \cos(\omega b) & \sin(\omega b) \\ \cosh(\omega b) & \sinh(\omega b) \end{bmatrix}.
\]
If $c_j(f)$ are such that
\[
    M\begin{bmatrix} c_1(f) \\ c_2(f)\end{bmatrix}=\begin{bmatrix}0 \\ -\int_0^1  [\Phi_1(1)\Phi_2(y)-\Phi_2(1)\Phi_1(y)]f(y)\,\mathrm{d}y \end{bmatrix},
\]
then, $u\in \operatorname{Dom}(D)$. We know that $\det(M)\not=0$, otherwise $i$ would be an eigenvalue of $D$ and this is impossible. Hence,
\[\begin{bmatrix} c_1(f) \\ c_2(f)\end{bmatrix}=M^{-1}\begin{bmatrix} 0 \\ \langle \Psi_2,f \rangle \end{bmatrix}\]
where $\Psi_2(y)=\Phi_2(1)\Phi_2(y)-\Phi_1(1)\Phi_2(y)$.

Now, in the expression \eqref{varpa}, the operators $f\longmapsto c_j(f)\Phi_j$ are rank one and the kernel
$[\Phi_1(x)\Phi_2(y)-\Phi_2(x)\Phi_1(y)]$ is bounded for all $(x,y)\in [0,1]^2$. Thus, indeed the integral operator given by the right hand side of  \eqref{varpa} is compact.
\end{proof}

\section*{Acknowledgements}
The work of LB was supported by COST Action CA18232. GF was supported by an EPSRC Research Associate grant. BP was partially supported by a Leverhulme Research Fellowship. DAS gratefully acknowledges support from the Quarterly Journal of Mechanics and Applied Mathematics Fund for Applied Mathematics.

\bibliographystyle{amsplain}
\bibliography{references}

\end{document}